\documentclass[12pt, a4paper,reqno]{amsart}
\usepackage{amsthm, amssymb, verbatim, cite, multirow}
\usepackage[utf8]{inputenc}
\usepackage[english]{babel}
\usepackage{hyperref}
\newtheorem{theorem}{Theorem}[section]
\newtheorem{cor}{Corollary}[section]
\newtheorem{lemma}{Lemma}[section]
\newtheorem{prop}{Proposition}[section]

\newtheorem{prob}{Problem}
\theoremstyle{definition}
\newtheorem{notation}{Notation}[section]
\newtheorem{remark}{Remark}[section]

\newcommand{\Aut}{\operatorname{Aut}}
\newcommand{\Out}{\operatorname{Out}}
\newcommand{\Inndiag}{\operatorname{Inndiag}}
\newcommand{\Outdiag}{\operatorname{Outdiag}}

\newcommand{\Soc}{\operatorname{Soc}}



\begin{document}

\title[On primes not adjacent to $2$ in the prime graphs]{\textsc{On primes not adjacent to $2$ in the prime graphs of almost simple groups}}

\author{M.A. Grechkoseeva}
\address{Maria A. Grechkoseeva, Sobolev Institute of Mathematics, Novosibirsk 630090, Russia}
\email{grechkoseeva@gmail.com}

\author{A.S. Kondrat'ev}

\address{Anatoly S. Kondrat'ev, N.N. Krasovskii Institute of Mathematics and Mechanics,  Yekaterinburg 620077, Box no. 82, Russia}
\email{a.s.kondratiev@imm.uran.ru}

\author{Nanying Yang}

\address{Nanying Yang, School of Mathematics and Data Science, Jiangnan University, Wuxi 214122, People’s Republic of China}
\email{yangny@jiangnan.edu.cn}

\begin{abstract}
The prime graph $GK(G)$ of a finite group $G$ is the graph whose vertices are primes dividing the order of $G$ and in which two distinct primes $r$ and $s$ are adjacent if and only if $G$ has an element of order $rs$. The $2$-independence number of $GK(G)$ is the maximal size of independent vertex subsets in $GK(G)$ containing the vertex~$2$. We determine all almost simple groups whose prime graph has a vertex $r$ not adjacent to $2$ and find the $2$-independence numbers of  all almost simple groups.


\smallskip

\end{abstract}

\maketitle

\section{Introduction}

Let $G$ be a finite group. Denote by $\pi(G)$ the set of all primes  dividing $|G|$ and by $\omega(G)$ the set of element orders of $G$. The \emph{prime graph}, or the Gruenberg--Kegel graph,  of $G$ is the graph $GK(G)$ with vertex set $\pi(G)$ in which two vertices $r$ and $s$ are adjacent if and only if $r\neq s$ and $rs\in\omega(G)$. Recall that a subset of vertices of a graph is a \emph{clique} if it consists of pairwise adjacent vertices, and a \emph{coclique}, or an independent subset, if it consists of pairwise non-adjacent vertices. Also recall that $G$ is an \emph{almost simple group} if $L\leq G\leq \Aut L$ for some non-abelian simple group $L$, in which case $L$ is the \emph{socle} of $G$.

It is clear that the structure of a finite group whose prime graph has non-adjacent vertices cannot be arbitrary, and of most importance are primes not adjacent to the vertex $2$. Following \cite{26YGKN},   we say that $G$ is a \emph{$W_r$-group} for an odd prime $r$ if $2,r\in\pi(G)$ but $2r\not\in\omega(G)$. 

Let $G$ be an arbitrary $W_r$-group. The main result of \cite{00CIY} states that $G$ has abelian Sylow $r$-subgroups. Also by \cite[Theorem 3]{00CIY}, if $G$ is not solvable, then $G$ has a~unique non-abelian composition factor; in other words, $G$ has a normal series $K\leq H\leq G$, where $K$ is the largest solvable normal subgroup of $G$, $S=H/K$ is a non-abelian simple group and $G/K$ is almost simple with socle $S$. The latter result was independently established in \cite{05Vas.t}, with the following important addition: $r\not\in\pi(G/H)$ and either $r\not\in\pi(K)$, or $S$ is specific, more precisely, $S\simeq A_7$ or $S\simeq PSL_2(q)$ for some odd $q$. Further improvement on the structure of non-solvable $W_r$-groups is obtained in \cite[Theorem 1.2]{26YGKN}, and this result largely reduces the problem of describing all non-solvable $W_r$-groups to 
\begin{enumerate}
\item determining almost simple $W_r$-groups, and 
\item investigating faithful 2\nobreakdash-\hspace{0pt}modular representations of quasisimple $W_r$-groups with fixed-point-free action of elements of order~$r$ (see \cite[Remark 1.3]{26YGKN}).
\end{enumerate}

In the present paper,  we deal with the first task.

\begin{prob}\label{p:asWr}
Determine almost simple $W_r$-groups for every odd prime $r$. 
\end{prob}

The structure of arbitrary $W_3$-groups is described in \cite{18KonMin.t}, in particular, the almost simple $W_3$-groups are listed in \cite[Corollary 1]{18KonMin.t}. The almost simple $W_5$-groups were found in   \cite[Theorem 1]{24GGKN.t}. If $r$ is arbitrary, then a reasonable restatement of Problem \ref{p:asWr} is as follows. 

\begin{prob}\label{p:asWr'}
Given an almost simple group $G$ and an odd prime $r$, determine whether $G$ is a $W_r$-group or not.
\end{prob}

The results of \cite{05VasVd.t} provide necessary and sufficient conditions of adjacency of vertices in the prime  graphs of simple groups (see also \cite{11VasVd.t} for corrections), so Problem \ref{p:asWr'} is solved for simple groups. 

We solve Problem \ref{p:asWr'} for all almost simple groups: Theorem \ref{t:1} covers almost simple groups of Lie type, while Propositions \ref{prop:alt} and \ref{prop:spor} are concerned with almost simple groups with alternating and sporadic socle, respectively. 
In the proof of Theorem \ref{t:1}, we use the result of \cite[Corollary 1]{18KonMin.t} to handle the case $r=3$, because $3$ is too specific to be considered alongside the others primes,  but we do not use the result of \cite[Theorem 1]{24GGKN.t} for the case $r=5$. Instead, we apply Theorem \ref{t:1} to describe almost simple $W_5$-groups of Lie type and compare this description with that of \cite[Theorem 1]{24GGKN.t} (see Corollary \ref{c:w5} and Remark \ref{rem:w5}).

Given a group $G$ of even order, denote by $\tau(G)$ the set of odd primes $r\in\pi(G)$ such that $2r\not\in\omega(G)$ (if $|G|$ is odd then $\tau(G)$ is not defined). It is clear that $G$ is a $W_r$-group if and only if $r\in\tau(G)$, so solving Problem \ref{p:asWr'}, we in fact deal with the sets $\tau(G)$ for almost simple groups $G$. For simple groups $L$, the sets $\tau(L)$ are found in \cite{05VasVd.t}, together with the numbers $t(2,L)$, where $t(2,G)$ denotes the maximal size of cocliques of $GK(G)$ containing the vertex~$2$ (and is referred to as the \emph{$2$-independence number} of $GK(G)$). Both the sets $\tau(L)$ and numbers $t(2,L)$, where  $L$ is simple, play a fundamental role in the problem of recognition of simple groups by the set of elements orders or the prime graph (see, for example, \cite[\S 2]{05Vas.t}), and one can expect that the sets $\tau(G)$ and numbers $t(2,G)$, where $G$ is almost simple, will be of equal importance for recognition of almost simple groups. This leads us to the following problem.

\begin{prob}\label{p:ast2}
Find $\tau(G)$ and $t(2,G)$ for all almost simple groups $G$.
\end{prob}

It is worth noting that for every group $G$ of even order, the subgraph induced by $\tau(G)$ in $GK(G)$ is a union of cliques \cite[Proposition 3.3]{24CGMY}, and $t(2,G)$ is just one more than the number of these cliques.
As we mentioned, Problem \ref{p:ast2} is solved for simple groups. Also $t(2,G)$ is known if the socle of $G$ is a Suzuki or Ree group \cite[Theorem 3]{23Gre.t}. 

We solve Problem \ref{p:ast2} for all almost simple groups: Theorem \ref{t:2} gives the answer for groups with Lie type socle,  and Propositions \ref{prop:alt} and \ref{prop:spor}  cover the case of  alternating or sporadic socle. 

\section{Preliminaries} \label{s:pre}

In this section, we collect general notation and necessary facts on primes not adjacent to~$2$, and also prove some auxiliary number-theoretical results.

\subsection{Group-theoretical notation}\label{ss:1}

Our group-theoretical notation and terminology are mostly standard and can be found in~\cite{08Isa}. In particular, $\Soc(G)$ denotes the socle of a group $G$, that is, the subgroup generated by minimal normal subgroups of $G$. If $\pi$ is a set of primes then $O_\pi(G)$ is the largest normal $\pi$-subgroup of $G$ and  $O^\pi(G)$ is the smallest normal subgroup of $G$ such that the corresponding quotient group is a $\pi$-group. As usual, $\pi'$ stands for the complement of $\pi$ in the set of all primes, and if $\pi=\{r\}$, then we write $O_{r}(G)$ and $O_{r'}(G)$ instead of $O_{\{r\}}(G)$ and $O_{\{r\}'}(G)$, respectively. Also, for brevity, if $G$ is an abelian group, then we write $(G)_{\pi}$ instead of $O_{\pi}(G)$.

For finite simple classical groups, we use the one-letter notation of the \textsc{Atlas} \cite{85Atlas}. In notation of classical groups, simple sporadic groups and simple groups of Lie type we follow Table 2, Table 1 and Table 5 of \cite{85Atlas}, respectively. Also we use the abbreviations $L_n^\varepsilon(q)$, $PGL_n^\varepsilon(q)$ and $E_6^\varepsilon(q)$, where $\varepsilon\in\{+,-\}$ and $L_n^+(q)=L_n(q)$, $L_n^-(q) = U_n(q)$, $PGL_ n^+(q) =PGL_n(q)$, $PGL_n^-(q) = PGU_n(q)$, $E_6^+(q)=E_6(q)$, and $E_6^-(q)={}^2E_6(q)$. 

Suppose that $L=L_2(q)$, where $q=p^f$, $p$ is odd and $f$ is even. Following \cite{26YGKN}, we write
$PGL_2^*(q)$ to denote the group $L\langle \delta\varphi\rangle$, where $L\langle \delta\rangle=PGL_2(q)$ and $\varphi$ is the field automorphism of $L$ of order $2$. 

\subsection{Number-theoretical notation and results} 
Given positive integers $a$ and $b$, we write $(a,b)$ for the greatest common divisor of $a$ and $b$. Also $\pi(a)$ denotes the set of prime divisors of $a$. If $\pi$ is a set of primes then $(a)_\pi$ is the greatest divisor $d$ of $a$ such that $\pi(d)\subseteq \pi(a)$. If $\pi=\{r\}$, then we write $(a)_r$ instead of $(a)_\pi$ and $(a)_{r'}$ instead of $(a)_{\pi'}$. We freely identify the signs $+$ and $-$ with the integers $1$ and $-1$, respectively, whenever they appear in arithmetic expressions. 

Let $a$ be an integer and $|a|>1$. An odd prime $r$ is called a \emph{primitive prime divisor} of $a^i-1$ if $r$ divides $a^i-1$ but does not divide $a^j-1$ for all $1\leq j<i$; and $2$ is a primitive divisor of
$a^i-1$ if $4$ divides $a^i-1$ but does not divide $a^j-1$ for all $1\leq j<i$. The set of all primitive divisors of $a^i-1$ is denoted by $R_i(a)$ and $r_i(a)$ stands for some element of $R_i(a)$ if it is not empty. It clear that $2$ is in $R_1(a)$ or $R_2(a)$ for odd $a$, 
depending on whether $a$ is congruent to $1$ or $-1$ modulo $4$.
The existence of primitive divisors for almost all pairs $i$ and $a$ was proved by Bang~\cite{86Bang} and Zsigmondy~\cite{Zs}.

\begin{lemma}[Bang--Zsigmondy]\label{l:bz}
Let $a$ be an integer with $|a| > 1$ and $i\in \mathbb N$.  If $R_i(a)= \varnothing$, then the pair $(i,a)$ is in $\{(1, 2),(1, 3),(6, 2)$, $(2, -2),(2, -3),(3, -2)\}$.
\end{lemma}

\begin{lemma}\label{l:bz0}
Let $p$ be a prime, $f,i\in\mathbb N$,  $\varepsilon\in \{+,-\}$, $q=p^f$, and suppose that $R_{i}(\varepsilon q)\neq \varnothing$. 
\begin{enumerate}
 \item If $R_{i}(\varepsilon q)\setminus\pi(f)=\varnothing$, then $q=8$ and $(i,\varepsilon)\in\{(2,+), (1,-)\}$.
 \item Let $f$ be even, $i\geq 3$ odd, $\varepsilon=+$ and set $q_0=p^{f/2}$.
 Then $R_{i}(q_0)\setminus\pi(f)\neq \varnothing$ and if $R_{2i}(q_0)\setminus\pi(f)=\varnothing$, then $q=4$ and $i=3$. 
 \end{enumerate}
\end{lemma}

\begin{proof}
(i) Suppose first that $\varepsilon=+$ and note that $R_{im}(p)\subseteq R_i(p^m)$. If $r\in R_{im}(p)$, then $r\equiv 1\pmod{im}$ and, in particular, $r\not\in\pi(m)$, so we may assume that  $R_{im}(p)=\varnothing$. Then $m>1$ and hence $im=6$ and $p=2$ by Lemma \ref{l:bz}. If $m=3$ and $i=2$, then $R_{i}(p^m)=\pi(m)$. If $m=2$ and $i=3$, or $m=6$ and $i=1$, then $7\in R_{i}(p^m)\setminus\pi(m)$.  

To handle the case $\varepsilon=-$, it suffices to observe that $R_i(-p^m)=R_i(p^m)$ if 
$i\equiv 0\pmod 4$, $R_i(-p^m)=R_{i/2}(p^m)$ if $i\equiv 2\pmod 4$, and $R_i(-p^m)=R_{2i}(p^m)$ if $i$ is odd.

(ii) Note that $R_{im/2}(p)\subseteq R_{i}(p^{m/2})$. Also if $r\in R_{im/2}(p)$, then $r>2$ and $r$ is coprime to $m$. If $R_{im/2}(p)=\varnothing$, then by (i) we have  $i=3$, $m=4$ and $p=2$, in which case $R_{i}(p^{m/2})=\{7\}$. Hence $R_{i}(q_0)\setminus\pi(f)\neq \varnothing$ in any case.

Since $R_{im}(p)\subseteq R_{2i}(p^{m/2})$ and $m>1$, we can argue as above and reduce to the case when $im=6$ and $p=2$. Then $i=3$, $m=2$ and $R_{2i}(p^{m/2})\subseteq \pi(m)$.
\end{proof}

For $\varepsilon \in\{+,-\}$, we set 
$$\Psi_4^\varepsilon(x) = x + \varepsilon\sqrt{2x} + 1,$$
$$ \Psi_6^\varepsilon(x) = x + \varepsilon\sqrt{3x} + 1,$$
$$\Psi_{12}^\varepsilon(x) = x^2 + \varepsilon\sqrt{2x^3} + x + \varepsilon\sqrt{2x} + 1.$$

\begin{lemma}\label{l:bzsr}
Let $p = 2$, $3$, or $2$, and $n=4$, $6$, or $12$, respectively. Let $f\geq 3$ be an odd integer and $\varepsilon\in\{+,-\}$. If $\pi(\Psi_n^\varepsilon(p^f))\cap R_{nf}(p)=\varnothing$, then $p=2$, $n=4$ and $f\in\{3,5\}$.
\end{lemma}

\begin{proof}
See \cite[Theorem 2]{23Gre.t}. 
\end{proof}

\begin{lemma}\label{l:bzsr0}
Let $p=2$, $3$, or $2$, and $n=4$, $6$, or $12$, respectively. Let $f\geq 3$ be an odd integer and $\varepsilon\in\{+,-\}$. If $\pi(\Psi_n^\varepsilon(p^f))\setminus \pi(f)=\varnothing$, then $p=2$, $n=4$, $m=5$ and $\varepsilon=-$.
\end{lemma}

\begin{proof}
If $r\in  \pi(\Psi_n^\varepsilon(p^f))\cap R_{nf}(p)$, then $r$ is coprime to $f$, so by Lemma \ref{l:bzsr}, we may assume that $p=2$, $n=4$ and $f\in\{3,5\}$. It remains to observe that $3\not\in \pi(2^3\pm 2^2+1)$, $\pi(2^5+2^3+1)=\{41\}$ and $\pi(2^5-2^3+1)=\{5\}$.
\end{proof}

\subsection{Prime graphs of simple groups} 

Given a finite group $G$ and $r,s\in\pi(G)$, we write $r\sim s$ if $r=s$ or $rs\in\omega(G)$. 
As we mentioned, for every $G$ of even order, the subgraph induced by $\tau(G)$ is a union of cliques, in other words, $\sim$ is an equivalence relation on $\tau(G)$. We will need this result only for simple groups.


\begin{lemma}\label{l:partition} Let $L$ be a finite simple group. Then $\sim$ is an equivalence relation on $\tau(L)$, $t(2,L)=1+|\tau(L)/{\sim}|$ and every coclique of $GK(L)$ containing $2$ can be included into a~coclique of size $t(2,L)$.  
\end{lemma}

\begin{proof}
 The first assertion follows from the adjacency criteria in the prime graphs of simple groups in \cite{05VasVd.t} (see, for example, \cite[Lemma 3.1]{24CGMY}). The remaining assertions are immediate corollaries of the first one.
\end{proof}

\begin{lemma}\label{l:simple} Suppose that $L$ is a non-abelian simple group,  $L\leq G\leq \Aut L$ and $r$ is an odd prime.
\begin{enumerate}
 \item If $G$ is a $W_r$-group, then  $r\not\in\pi(G/L)$ and, in particular, $L$ is a $W_r$-group.
 \item If $L$ is a $W_r$-group, $r\not\in\pi(G/L)$ and $G$ is not a $W_r$-group, then   there is an involution $t\in G\setminus L$ such that $r\in\pi(C_L(t))$.
\end{enumerate}

\end{lemma}
\begin{proof}
Item (i) follows, for example, from \cite[Lemma 1.2]{05Vas.t}. To prove (ii), 
let $g\in G\setminus L$ be an element of order $2r$. Since $r\not\in\pi(G/L)$, it follows that $g^2\in L$,  and so $t=g^r\not\in L$.  
\end{proof}





\begin{lemma}\label{l:w3}
Let $G$ be an almost simple group. Then $G$ is a $W_3$-group if and only if $O^{\{2,3\}'}(G)$ is one of the following groups: $L_2(2^f)$, $L_2(3^f)$, $PGL_2(3^f)$, $PGL_2^*(3^f)$ ($f$ is even), $L_2(q)$, where $q\equiv 
\pm 5 \pmod {12}$, $L_3^\varepsilon (2^f)$, where $(2^f-\varepsilon)_3 \leq 3$. 
\end{lemma}

\begin{proof}
See \cite[Theorem 2]{18KonMin.t}. 
\end{proof}

\section{Almost simple $W_r$-groups of Lie type}\label{s:wrlie}

In this section, we introduce notation for some outer automorphisms of simple groups of Lie type, mostly for involutions, and then prove Theorem \ref{t:1}.

Throughout the section $L$ is a simple group of Lie type and $\overline{\phantom{x}}:\Aut L\rightarrow \Out L$ is the natural homomorphism. Recall that we denote simple groups of Lie type according to \cite[Table 5]{85Atlas}, so $L={}^d\Sigma(q)$, where $\Sigma$ is a root system, $d$ is empty, $2$, or $3$, and $q=p^f$. Note that this notation slightly differs from  the notation of Theorem 2.5.12 of  \cite{98GorLySol}, which we will cite in the next paragraph: the simple Suzuki and Ree groups are denoted by $^2\Sigma(p^f)$ in \cite{85Atlas}  and by $^2\Sigma(p^{f/2})$ in the theorem.

Assume that $L\neq{}^2F_4(2)'$. By \cite[Theorem 2.5.12]{98GorLySol}, we have $\Aut L = \Inndiag L\rtimes \Phi\Gamma$, where $\Inndiag L$, $\Phi$, and $\Gamma$ are respectively the groups of inner-diagonal, standard field, and standard graph automorphisms of $L$. Furthermore, $\Phi$ is a cyclic group of order $f$ if $d=2$, $\Sigma=B_2$, $G_2$, $F_4$ and $p=2,3,2$, respectively, and  $\Phi$ is a cyclic group of order $df$ otherwise. Also $\Phi\Gamma=\Phi\times \Gamma$ if $d=1$ and $\Sigma$ is simply laced; $\Phi\Gamma$ is cyclic with $|\Phi\Gamma:\Phi|=2$ if $d=1$, $\Sigma=B_2$, $G_2$, $F_4$ and $p=2,3,2$, respectively; and $\Gamma=1$ otherwise.
As usual, $\Outdiag L=\Inndiag L/L$. 

If $L={}^2F_4(2)'$, then we set $\Outdiag L=\overline \Phi=1$ and $\overline{\Phi\Gamma}=\Out L$. 

In what follows, $G$ is an almost simple group with socle $L$. 
 By Lemma \ref{l:simple}, if $G$ is a $W_r$-group,  then $L$ is a $W_r$-group and $r\not\in\pi(\overline G)$. We can say more about $r$ and $G$. 

\begin{lemma}\label{l:general} Let $G$ be a $W_r$-group for some prime $r$.  The following hold.
\begin{enumerate}
 \item $r\not\in\pi(\Outdiag L)$ unless $r=3$, $L=L_3^\varepsilon(q)$, $q$ is even and $(q-\varepsilon)_3=3$.
  \item Sylow $2$-subgroups of $\overline G$ are cyclic. 
  \item If $L\neq L_2(q)$, then $|\overline G\cap\Outdiag L|$ is odd.
\end{enumerate}
 \end{lemma}
 
\begin{proof} 
(i) If $r$ divides $|\Outdiag L|$, then $L=L_n^\varepsilon(q)$ since $6\in \omega(E_6^\pm(q))$. By \cite[Tables 4 and 6]{05VasVd.t}, a prime divisor $r$ of $(n,q-\varepsilon)$ is not adjacent to $2$ in $GK(L)$ if and only if $n=3$, $q$ is even, $(q-\varepsilon)_3=3$ and $r=3$.

(ii)  If $r=3$, this follows from Lemma \ref{l:w3}, so we may assume that $r>3$. 

Suppose that $L\neq L_2(q)$ or $r\neq p$. Then a Sylow $r$-subgroup $R$ of $L$ is cyclic (see, for example, \cite[Lemma 2.3]{22GreVas}). Hence $N_G(R)/C_G(R)$ is also cyclic. Since $|C_G(R)|$ is odd and $N_G(R)$ covers $\overline G$, the claim follows in this case.

Let $L=L_2(q)$ and $r=p$. Every field automorphism of $L$ centralizes an element of order~$p$, so $|\overline G\cap\overline \Phi|$ is odd. Thus the order of a Sylow 2-subgroup of $\overline G$ is at most 2. 

(iii) See \cite[Theorem 2(i)]{22GreVas}.
\end{proof}

The next easy corollary of Lemmas \ref{l:simple} and \ref{l:general} 
provides a necessary and sufficient condition for $G$ to be a $W_r$-group.

\begin{lemma} \label{l:cr} Let $r$ be an odd prime. Then $G$ is a $W_r$-group if and only if $r\not\in\pi(\overline G)$,  Sylow $2$-subgroups of $\overline G$ are cyclic and one of the following holds:
\begin{enumerate}
 \item $|\overline G|$ is odd and $L$ is a $W_r$-group;
 \item $|\overline G|$ is even and $G_0$ is a $W_r$-group, where $L\leq G_0\leq G$ and $|G_0:L|=2$. 
\end{enumerate}
\end{lemma}

By Lemma \ref{l:cr}, to determine all almost simple $W_r$-groups with socle $L$, we need to find all almost simple $W_r$-group $G_0$ such that $|G_0:L|=2$. 
Let $\alpha$ be an involution in $\Out L$ and take some preimage $\dot{\alpha}$ of $\alpha$ in $\Aut L$. Let $\mathcal C_\alpha$ be the set of representatives of the $\Aut L$-conjugacy classes of involutions in $L\dot\alpha$ and let $\pi_\alpha$ be the union of $\pi(C_L(t))$ over $t\in \mathcal C_\alpha$. By Lemma \ref{l:simple}(ii), we see that $G_\alpha=\langle L,\dot\alpha\rangle$ is a $W_r$-group if and only if $L$ is a $W_r$-group and $r\not\in \pi_\alpha$. 

Before proceed with calculation of the sets $\pi_\alpha$, we introduce notation for some elements of $\Out L$ and list the conjugacy classes of involutions of $\Out L$ outside $\Outdiag L\rtimes \overline \Gamma$ in the cases when these classes are not clear.

\begin{notation}\label{n:auto}

~

\begin{enumerate}
 \item $\varphi$ is the only involution of $\overline \Phi$ if $d\neq 2$ and $q$ is a square (so writing $\varphi$ we always assume that $f$ is even and put $q_0=p^{f/2}$);
 \item $\psi$ is the only involution of $\overline \Phi\overline \Gamma$ if $d=1$, $\Sigma=B_2$, $G_2$, $F_4$ and $p=2,3,2$, respectively, and $f$ is odd;
  \item $\gamma$ is the image of the standard graph automorphism of order $2$
  if $d=1$, or the only involution of $\overline \Phi$ if $d=2$ ($\Sigma$ is simply laced in both cases); 
 \item $\gamma'$ is an involution in $\Outdiag L\gamma$ not conjugate to $\gamma$ (only if $q$ is odd and $L$ is either $L_n^\pm(q)$, where $n$ is even, or $O_{2n}^\varepsilon(q)$, where $n$ is odd, or $O_{2n}^-(q)$, where $n$ is even);
 \item $\delta$ is the only involution of $\Outdiag L$ if $\Outdiag L$ is cyclic of even order;
 \item $\tau$ is a generator of $O_3(\overline \Gamma)$ when $L=O_8^+(q)$;
 \item $\upsilon$ is a generator of $\Out L$ when $L={}^2F_4(2)'$.
\end{enumerate}
\end{notation}

Observe that $\overline{PGL_2^*(q)}=\langle \delta\varphi\rangle$ in this notation.

\begin{lemma}\label{l:involutions}
Let $L=L_n(q)$, where $n$ is even, $q$ is odd, $f$ is even, and let $\alpha$ be involution in $\Out L\setminus (\Outdiag L\rtimes \overline \Gamma)$. If $(n)_2\geq (q-1)_2$, then $\alpha$ is conjugate to $\varphi$ or $\varphi\gamma$.  If $(n)_2<(q-1)_2$, then $\alpha$ is conjugate to  $\varphi$, $\varphi\gamma$,
$\varphi\delta$, or $\varphi\gamma\delta$.
\end{lemma}

\begin{proof}
Let $(\Outdiag L)_2=\langle \eta\rangle$, $\varphi_+=\varphi$ and $\varphi_-=\varphi\gamma$.  Then $\eta^{\varphi_\epsilon}=\eta^{\epsilon q_0}$, where $q_0=p^{f/2}$. 

The inequality $(n)_2\geq (q-1)_2$ yields $(n,q-1)_2=(q-1)_2$. So $\varphi_\epsilon\eta^i$ is an involution if and only if $i(\epsilon q_0+1)$ is divisible by $(q-1)_2$ or, equivalently, $i$ is divisible by $(\epsilon q_0-1)_2$. On the other hand, $(\varphi_\epsilon)^\eta=\varphi_\epsilon\eta^{1-\epsilon q_0}$. Since $(1-\epsilon q_0, (q-1)_2)=(q_0-\epsilon)_2$, the group $\langle \eta\rangle$ acts transitively on the set of involutions of the form $\varphi_\epsilon \eta^i$.

Let $(n)_2<(q-1)_2$ and choose $\epsilon$ so that $q_0\equiv \epsilon\pmod 4$. Then $\varphi_\epsilon$ centralizes $\eta$ and hence $\alpha$ is conjugate to an element of $(\langle \eta\rangle\rtimes \langle \gamma\rangle)\times \langle \varphi_\epsilon\rangle$. It follows that $\alpha$ is conjugate to one of $\varphi_\epsilon$, $\varphi_\epsilon\gamma$, $\varphi_\epsilon\delta$, and $\varphi_\epsilon\gamma\delta$, as stated. 
\end{proof}

Now we return to involutions $\alpha$ such that $G_\alpha$ is a $W_r$-group. 
By Lemma \ref{l:general}(iii), if $L\neq L_2(q)$, then we may assume that $\alpha\not\in\Outdiag L$. By \cite[Proposition 4.9.1(d)]{98GorLySol}, if $\alpha\in \Outdiag L\beta$, where $\beta=\varphi$ or $\gamma\varphi$, and $\alpha$ is not $\Out L$-conjugate to $\beta$, then there are no involutions in $G_\alpha\setminus G$, and so $\pi_\alpha=\varnothing$. Hence we are left with the cases when $\alpha\in\overline\Phi\overline\Gamma$ or $\alpha=\gamma'$.
In the next lemma we collect information about $\pi_\alpha$ for such $\alpha$.

\begin{lemma} \label{l:pia}
Let  $\alpha\in \Out L$ and $|\alpha|=2$.
\begin{enumerate}
 \item If $\alpha=\varphi$, then $\pi_\alpha=\pi(\Sigma(q_0))$. 
 
 \item If $\alpha=\varphi\gamma$, then $\pi_\alpha=\pi({}^2\Sigma(q_0))$. 
 
 \item If $\alpha=\psi$, then $\pi_\alpha=\pi({}^2\Sigma(q))$.
 
 \item If $\alpha=\gamma$ and $L=L^\pm_n(q)$, then $\pi_\alpha=\pi(S_n(q))$ or $\pi(O_n(q))$, depending on whether $n$ is even or odd.
 
 \item If $\alpha=\gamma$ and $L=O^\pm_{2n}(q)$, then $\pi_\alpha=\pi(S_{2n-2}(q))$.
 
 \item If $\alpha=\gamma$ and $L=E_6^\pm(q)$, then $\pi_\alpha=\pi(F_4(q))$.
 
 \item If $\alpha=\gamma'$ and $L=L_n^\varepsilon(q)$, with $q$ odd and $n$ even, then $\pi_\alpha=\pi(O^\epsilon_n(q))$, where $\epsilon=-$ if $n\equiv 0\pmod 4$ or $q\equiv \varepsilon\pmod 4$ and $\epsilon=+$ otherwise.
 
 \item If $\alpha=\gamma'$ and $L=O_{2n}^\varepsilon(q)$, with $q$ and $n$ odd, then $\pi_\alpha=\pi(S_{n-1}(q^2))$.
 \item If $\alpha=\gamma'$ and $L=O_{2n}^-(q)$, with $q$ odd and $n$ even, then $\pi_\alpha=\varnothing$.
 \item If $\alpha=\upsilon$, then $\pi_\alpha=\varnothing$. 
\end{enumerate}
 
\end{lemma}

\begin{proof}

The claim follows from \cite[Proposition 4.9.1]{98GorLySol} in (i)--(iii); \cite[Table 4.5.1]{98GorLySol} in other cases if $q$ is odd;  \cite[(19.8)]{76AscSei} in (iv) and (iv) if $q$ is even;  \cite[(8.7), (8.10)]{76AscSei} in (v) if $q$ is even; and \cite{85Atlas} in (x).

As an example of how to find $\pi_\alpha$ using these results, we consider  Item (iv). Suppose that $q$ is odd. If $n$ is odd, then by \cite[Table 4.5.1]{98GorLySol}, there is the only conjugacy class of involutions in $L\dot\gamma$, with centralizer of type $B_{(n-1)/2}(q)$. If $n$ is even, then there are two conjugacy classes of involutions with centralizes of type $C_{n/2}(q)$ and $D_{n/2}^\epsilon(q)$ for some $\epsilon$. It remains to note that $\pi(D_{n/2}^\epsilon(q))\subseteq\pi(C_{n/2}(q))$. Similarly, if $q$ is even, then by \cite[(19.8)]{76AscSei}, we have $\pi_\gamma=\pi(SO_n(q))$ or $\pi(Sp_n(q))$, depending on whether $n$ is odd or even.
\end{proof}

In the proof of the next proposition, we will work with the sets $\tau(L)$. The description of these sets are given in \cite[Tables 4 and 5]{05VasVd.t} for even $q$ and \cite[Tables 6 and 7]{05VasVd.t} for odd $q$, and we use these results without further references. Note that there is a misprint in Table 4 of \cite{05VasVd.t} for $L={}^2A_{n-1}(q)$: if $2=(n)_2=(q+1)_2$, then it should be $\{2,r_{2n-2}(q)\}$ in the last column, not $\{2,r_{n/2}(q)\}$. 

\begin{table}[t]
\caption{Almost simple $W_r$-groups $G$ of Lie type in characteristic $2$ with $|G/\Soc(G)|=2$}\label{tab:2even}

\renewcommand{\arraystretch}{1.3}
\begin{tabular}{|l|l|l|l|l|}
\hline 
$L=\Soc(G)$& \multicolumn{2}{l|}{Conditions on $L$}&$r$&$\alpha$\\
\hline 
$L_n(q)$& $n=2$ & $f$ even & $R_2(q)$ &$\varphi$\\
\cline{2-5}
& $n\geq 3$ odd & & $R_{n}(q)$&$\gamma$\\
\cline{3-5}
&& $f$ even& $R_{n-1}(q)\cup R_{2n}(q_0)$&$\varphi$\\ 
&  & & $R_{n-1}(q)\cup R_n(q_0)$&$\varphi\gamma$\\
\cline{2-5}
& $n\geq 4$ even & &$R_{n-1}(q)$&$\gamma$\\
\cline{3-5}
&& $f$ even&$R_n(q)\cup R_{2n-2}(q_0)$&$\varphi$\\ 
&  & &$R_n(q)\cup R_{n-1}(q_0)$&$\varphi\gamma$\\

\hline
$S_{2n}(q)$&$n\geq 2$& $f$ even &$R_{2n}(q)$&$\varphi$\\ 
\hline
$O_{2n}^+(q)$&$n\geq 4$ even &$f$ even&$R_{2n-2}(q)$&$\varphi$, $\varphi\gamma$\\ 
\cline{2-5}
&$n\geq 5$ odd &   &$R_{n}(q)$&$\gamma$\\
\cline{3-5}
&&$f$ even &$R_{2n-2}(q)\cup R_{2n}(q_0)$&$\varphi$\\ 
& & &$R_{2n-2}(q)\cup R_{n}(q_0)$&$\varphi\gamma$\\ 

\hline
$U_n(q)$& \multicolumn{2}{l|}{$n\geq 3$ odd}& $R_n(-q)$&$\gamma$\\
\cline{2-5}
&\multicolumn{2}{l|}{$n\geq 4$ even}& $R_{n-1}(-q)$&$\gamma$\\
\hline
$O_{2n}^-(q)$&\multicolumn{2}{l|}{$n\geq 4$}  &$R_{2n}(q)$&$\gamma$\\ 
\hline 
$G_2(q)$&\multicolumn{2}{l|}{$f$ even}& $R_{6}(q)$&$\varphi$\\
\hline

$^3D_4(q)$&\multicolumn{2}{l|}{$f$ even}&$R_{12}(q)$&$\varphi$\\
\hline
$F_4(q)$&\multicolumn{2}{l|}{$f$ even}&$R_{12}(q)\cup R_8(q)$&$\varphi$\\
\cline{2-5}
&\multicolumn{2}{l|}{$f$ odd}&$R_8(q)$&$\psi$\\
\hline
$^2F_4(2)'$&\multicolumn{2}{l|}{} &$13$&$\upsilon$\\
\hline
$E_6(q)$&  \multicolumn{2}{l|}{} &$R_{9}(q)$&$\gamma$\\
\cline{2-5}
& \multicolumn{2}{l|}{$f$ even}&$R_{12}(q)\cup R_8(q)\cup R_{18}(q_0)$&$\varphi$\\ 
&  \multicolumn{2}{l|}{} &$R_{12}(q)\cup R_8(q)\cup R_{9}(q_0)$&$\varphi\gamma$\\

\hline
$^2E_6(q)$& \multicolumn{2}{l|}{} & $R_{18}(q)$&$\gamma$\\
\hline
$E_7(q)$&\multicolumn{2}{l|}{$f$ even}& $R_{18}(q)\cup R_{14}(q)$&$\varphi$\\
\hline
$E_8(q)$&\multicolumn{2}{l|}{$f$ even}& $R_{20}(q)\cup R_{24}(q)\cup R_{30}(q)$&$\varphi$\\
\hline
\end{tabular}
\end{table}

\begin{prop}\label{p:2even}  Let $p=2$, $r\geq 3$ and $\overline G=\langle \alpha\rangle$, where $\alpha\in \Out L$ is an involution. Then $G$ is a $W_r$-group if and only if, up to conjugacy, $\alpha$ is as in Table {\rm \ref{tab:2even}}. 
 \end{prop}
 
\begin{proof}

As we explained above, $G$ is a $W_r$-group if and only if $r\in \tau(L)\setminus \pi_\alpha$. So we need to calculate the set $\tau(L)\setminus \pi_\alpha$. Since $|\Outdiag L|$ is odd, every involution is conjugate to an element of $\overline \Phi\overline \Gamma$.

We give a detailed proof for the cases when $L=L_n^\varepsilon(q)$ or $F_4(q)$, the other cases are similar. Observe that for $q=q_0^2$, we have  $R_n(q)=R_n(q_0)\cup R_n(-q_0)=R_n(q_0)\cup R_{2n}(q_0)$ if $n\geq 3$ is odd and $R_n(q)=R_{2n}(q_0)$ if $n\geq 2$ is even. 

\textbf{1.} Let $L=L_2(q)$, where $q\geq 4$. Then $\tau(L)=R_1(q)\cup R_2(q)$ and $\alpha=\varphi$. By Lemma \ref{l:pia}(i), we have $\pi_{\varphi}=\pi(L_2(q_0))$. Since $q_0^2-1=q-1$, it follows that $\tau(L)\setminus \pi_\varphi=R_2(q)$. 

Let $L=L_n^\varepsilon(q)$, where $n\geq 3$. By Lemma \ref{l:w3}, there is no $W_3$-groups $G$ such that $|G/L|=2$. It is clear that $3\not\in R_i(q)$ if $i>2$ and $3\not\in R_2(q)$ if $f$ is even, so $3$ does not appear in Table~\ref{tab:2even}, so for $r=3$ the claim follows. 

Let $r>3$. We have $\tau(L)\setminus\pi(3,q-\varepsilon)=R_n(\varepsilon q)\cup R_{n-1}(\varepsilon q)$. 

Suppose that $\varepsilon=+$. Then $\alpha$ is one of $\varphi$, $\varphi\gamma$, $\gamma$ and we apply (i), (ii) or (iv) of Lemma~\ref{l:pia}, respectively, to obtain $\pi_\alpha$. 

If $n$ is odd, then $$R_n(q)\cap \pi(L^\epsilon_n(q_0))= R_n(\epsilon q_0)\text{ and } R_n(q)\cap \pi(O_n(q))=\varnothing,$$ $$R_{n-1}(q)\cap \pi(L^\pm_n(q_0))= \varnothing \text{ and } R_{n-1}(q)\subseteq  \pi(O_n(q)).$$

If $n$ is even, then $$R_{n}(q)\cap \pi(L^\pm_n(q_0))= \varnothing \text{ and }R_{n}(q)\subseteq  \pi(S_n(q)),$$ $$R_{n-1}(q)\cap \pi(L_n^\epsilon(q_0))= R_{n-1}(\epsilon q_0)\text{ and }R_{n-1}(q)\cap \pi(S_n(q))=\varnothing.$$ 

Suppose that $\varepsilon=-$. Then $\alpha=\gamma$. If $n$ is odd, then $$R_n(-q)\cap \pi(O_n(q))=R_{2n}(q)\cap \pi(O_n(q))=\varnothing, \text{ and }R_{n-1}(-q)\subseteq \pi(O_n(q)).$$ 

If $n$ is even then $$R_n(-q)\subseteq  \pi(S_n(q))\text{ and }R_{n-1}(-q)\cap\pi(S_n(q))=\varnothing.$$

\textbf{2.} Let $L=F_4(q)$, where $q=2^f\geq 4$. Then $\tau(L)=R_8(q)\cup R_{12}(q)$. If $f$ is even, then the only involution in $\overline{\Phi\Gamma}$ is $\varphi$ and $\pi_{\varphi}=\pi(F_4(q_0))$. If $f$ is odd, then the only involution is $\psi$ and $\pi_{\psi}=\pi({}^2F_4(q))$. In either case, $R_8(q)\cap \pi_\alpha=\varnothing$. Since $R_{12}(q)\cap \pi(F_4(q_0))=\varnothing$ and $R_{12}(q)\subseteq \pi({}^2F_4(q))$, the claim follows. 
\end{proof}

\begin{table}
 
 \caption{Almost simple $W_r$-groups $G$ of Lie type in odd characteristic $p$ with $|G/\Soc(G)|=2$}\label{tab:2odd}

\renewcommand{\arraystretch}{1.3}
\begin{tabular}{|l|l|l|l|l|}
\hline 
$L=\Soc(G)$& \multicolumn{2}{l|}{Conditions on $L$}&$r$&$\alpha$\\
\hline 
$L_n(q)$& $n=2$& & $p$& $\delta$\\ 
\cline{3-5}
&&$f$ even &$p$ &$\delta\varphi$\\
&  && $R_2(q)$ &$\varphi$, $\delta\varphi$\\
\cline{2-5}
& $n\geq 3$ odd &  & $R_{n}(q)$&$\gamma$\\
\cline{3-5}
&&$f$ even& $R_{2n}(q_0)$&$\varphi$\\ 
& & & $R_n(q_0)$&$\varphi\gamma$\\
\cline{2-5}
&$n\geq 4$ even& $(n)_2\geq (q-1)_2$& $R_{n-1}(q)$& $\gamma$, $\gamma'$\\
\cline{2-5}
& $n\geq 4$ even,&$(n)_2<(q-1)_2$ & $R_n(q)$&$\alpha\not\in \Outdiag L\rtimes \overline \Gamma$\\ 
\cline{3-5}
&$f$ even&$(n)_2=(q-1)_2$&$R_n(q)\cup R_{2n-2}(q_0)$&$\varphi$\\ 
&&  & $R_n(q)\cup R_{n-1}(q_0)$&$\varphi\gamma$\\
\cline{3-5}
&& $(n)_2>(q-1)_2$  &$R_{2n-2}(q_0)$&$\varphi$\\
&&&$R_{n-1}(q_0)$&$\varphi\gamma$\\
\hline
$U_n(q)$& \multicolumn{2}{l|}{$n\geq 3$ odd} &$R_n(-q)$&$\gamma$\\
\cline{2-5}
& $n\geq 4$ even& $(n)_2\geq (q+1)_2$ & $R_{n-1}(-q)$&$\gamma$, $\gamma'$\\
\cline{3-5}
&& $2=(n)_2<(q+1)_2$ & $R_{n}(-q)$&$\gamma'$\\
\hline
$S_{2n}(q)$&\multicolumn{2}{l|}{$n\geq 2$, $f$ even} & $R_{2n}(q)$&$\varphi$, $\delta\varphi$\\
\hline
$O_{2n+1}(q)$&\multicolumn{2}{l|}{$n\geq 3$, $f$ even} & $R_{2n}(q)$&$\varphi$, $\delta\varphi$\\
\hline
$O_{2n}^+(q)$&\multicolumn{2}{l|}{$n\geq 4$, $f$ even}  &$R_{2n-2}(q)$ &$\alpha\not\in \Outdiag L\rtimes \overline \Gamma$\\ 
\cline{2-5}
&\multicolumn{2}{l|}{$n\geq 5$ odd, $q\not\equiv 1\pmod 8$} &$R_{n}(q)$&$\gamma$, $\gamma'$\\
\hline
$O_{2n}^-(q)$&\multicolumn{2}{l|}{$n\geq 4$ even or $q\not\equiv 7\pmod 8$} &$R_{2n}(q)$&$\gamma$, $\gamma'$\\ 
\hline 
$G_2(q)$&\multicolumn{2}{l|}{$f$ even}& $R_{6}(q)$&$\varphi$\\
\cline{2-5}
&\multicolumn{2}{l|}{$p=3$, $f$ odd} & $R_3(q)$&$\psi$\\
\hline
$F_4(q)$&\multicolumn{2}{l|}{$f$ even}&$R_{12}(q)$&$\varphi$\\
\hline
$^3D_4(q)$&\multicolumn{2}{l|}{$f$ even}&$R_{12}(q)$&$\varphi$\\
\hline
$E_6(q)$&\multicolumn{2}{l|}{}&$R_{9}(q)$&$\gamma$\\
\cline{2-5}
&\multicolumn{2}{l|}{$f$ even} &$R_{12}(q)\cup R_{18}(q_0)$&$\varphi$\\ 
& \multicolumn{2}{l|}{} & $R_{12}(q)\cup R_{9}(q_0)$&$\varphi\gamma$\\

\hline
$^2E_6(q)$& \multicolumn{2}{l|}{} & $R_{18}(q)$&$\gamma$\\
\hline
$E_7(q)$&\multicolumn{2}{l|}{$f$ even}& $R_{18}(q)\cup R_{14}(q)$&$\varphi$, $\delta\varphi$\\
\hline
$E_8(q)$&\multicolumn{2}{l|}{$f$ even}& $R_{20}(q)\cup R_{24}(q)\cup R_{30}(q)$&$\varphi$\\
\hline
\end{tabular}
\end{table}

\begin{prop}  \label{p:2odd} Let $p\neq 2$, $r\geq 3$ and $\overline G=\langle \alpha\rangle$, where $\alpha\in \Out L$ is an involution. Then $G$ is a $W_r$-group if and only if, up to conjugacy, $\alpha$ is as in Table {\rm \ref{tab:2odd}}. 
 \end{prop}
 
\begin{proof}
If $\alpha$ is conjugate to an element of $\overline \Phi\overline \Gamma$, then we take $\pi_\alpha$ from Lemma \ref{l:pia} and just calculate $\tau(L)\setminus \pi_\alpha$ as in the proof of Proposition \ref{p:2even}. So we may assume that $|\Outdiag L|$ is even. 

\textbf{1.} Let $L$ be one of $S_{2n}(q)$, $O_{2n+1}(q)$ or $E_7(q)$. Then $\alpha\in \langle \delta\rangle \times (\overline \Phi)_2$. As we mentioned, $\alpha\neq \delta$, so we may assume that $f$ is even. Then $q\equiv 1\pmod 4$, and hence $\tau(L)=R_{2n}(q)$ or $R_{18}(q)\cup R_{14}(q)$. If $\alpha=\varphi$, then $\pi_\alpha=\pi(\Sigma(q_0))$. If $\alpha=\delta\varphi$, then $\pi_\alpha=\varnothing$. In either case,  $\tau(L)\cap \pi_\alpha=\varnothing$.

\textbf{2.} Let $L=O_{2n}^-(q)$. Then $\alpha=\gamma$ or $\gamma'$.
Also $\tau(L)$ is equal to $R_{2n}(q)$ if $n$ is even or $q\equiv 1\pmod 4$, $R_{2n}(q)\cup R_{2n-2}(q)$ if $n$ is odd and $q\equiv 3\pmod 8$, 
and $R_{2n-2}(q)$ if $n$ is odd and $q\equiv 7\pmod 8$. Applying Lemma \ref{l:pia}, we see that $R_{2n}(q)\cap \pi_\alpha=\varnothing$, $R_{2n-2}(q)\subseteq\pi_\gamma$, and if $n$ is odd, then $R_{2n-2}(q)\subseteq\pi_{\gamma'}$.

\textbf{3.} Let $L=U_n(q)$, where $n$ is even. Then $\alpha=\gamma$ or $\gamma'$. We have $\pi_{\gamma'}=\pi(O_{n}^+(q))$ if $(n)_2=(q+1)_2=2$, and $\pi_{\gamma'}=\pi(O_{n}^-(q))$ otherwise. Also  $\pi_\gamma=\pi(S_n(q))$.

Suppose that $(n)_2<(q+1)_2$. Then $\tau(L)=R_n(-q)$, and   hence
$\tau(L)\subseteq \pi_\gamma$. If $(n)_2=2$, then $R_n(-q)=R_{n/2}(q)$, and so $\tau(L)\cap \pi_{\gamma'}=\varnothing$. If $(n)_2>2$, then  $R_n(-q)=R_{n}(q)$ and $\tau(L)\subseteq \pi_{\gamma'}$. 

Suppose that $(n)_2>(q+1)_2$ or $2=(n)_2=(q+1)_2$. Then $\tau(L)=R_{n-1}(-q)=R_{2n-2}(q)$. We see that $\tau(L)\cap \pi_\alpha=\varnothing$. 

If $2<(n)_2=(q+1)_2$, then $\tau(L)=R_n(-q)\cup R_{n-1}(-q)=R_n(q)\cup R_{2n-2}(q)$, and hence $\tau(L)\setminus \pi_\alpha=R_{n-1}(-q)$.

\textbf{4.} Let $L=O_{2n}^+(q)$, where $n$ is even and $q\equiv \epsilon\pmod 4$. Then $\tau(L)=R_{2n-2}(\epsilon q)$ and so $\tau(L)\subseteq \pi_\gamma$. This yields $\alpha\not\in\Outdiag L\rtimes \overline \Gamma$. 
Hence $f$ is even and $\epsilon=1$. Since $\overline \Phi$ centralizes $\Outdiag L\rtimes \overline \Gamma$, it follows that $\alpha$ is conjugate to a product of $\varphi$ and some  $\beta\in\Outdiag L\rtimes \overline \Gamma$ such that $\beta^2=1$.
If $\alpha$ is conjugate to $\varphi$ or $\varphi\gamma$, then $\pi_\alpha\cap \tau(L)=\varnothing$. Otherwise $\pi_\alpha=\varnothing$.  

\textbf{5.} Let $L=O_{2n}^+(q)$, where $n$ is odd. 
Suppose that $q\equiv 1\pmod 8$. Then $\tau(L)=R_{2n-2}(q)$. In particular, $\tau(L)\subseteq \pi_\alpha$ for $\alpha=\gamma,\gamma'$. So $\alpha\not\in\Outdiag L\rtimes \overline \Gamma$. Arguing as in the previous case, we conclude that $f$ is even and $R_{2n-2}(q)\cap\pi_\alpha=\varnothing$ for all such $\alpha$.

If $q\not\equiv 1\pmod 8$, then $R_n(q)\subseteq \tau(L)\subseteq R_n(q)\cup R_{2n-2}(q)$ and $\alpha$ is one of $\gamma$ or $\gamma'$. In either case, $\tau(L)\setminus\pi_\alpha=R_n(q)$. 

\textbf{6.} Let $L=L_n(q)$, where $n\geq 4$ even. If $\alpha$ is $\gamma$ or $\gamma'$, then $\tau(L)\setminus \pi_\alpha$ is equal to $R_{n-1}(q)$ if $(n)_2\geq (q-1)_2$ and empty otherwise. Suppose that $\alpha\not\in\Outdiag L\rtimes \overline \Gamma$.
Then $f$ is even and, in particular, $(q-1)_2>4$. 

If $(n)_2<(q-1)_2$, then $\tau(L)=R_n(q)$. By Lemma \ref{l:involutions}, we may assume that $\alpha$ is one of $\varphi$, $\varphi\gamma$, $\varphi\delta$, $\varphi\gamma\delta$. 
In any case, $R_n(q)\cap\pi_\alpha=\varnothing$. 

If $(n)_2\geq (q-1)_2$, then every involution outside $\Outdiag L\rtimes \overline \Gamma$ is conjugate to $\varphi$ or $\gamma\varphi$. Calculating $\tau(L)\setminus \pi(L_n^\pm(q_0))$ completes the proof in this case. 

\textbf{7.} Let $L=L_2(q)$. If $\alpha=\delta$, then $\tau(G)=\{p\}$ since $\omega(PGL_2(q))$ consists of the divisors of $p$, $q-1$, and $q+1$. So we may assume that $f$ is even, and hence $\tau(L)=R_2(q)\cup\{p\}$, and $\alpha=\delta\varphi$ or $\varphi$. We have $\pi_\alpha=\varnothing$ in the former case and $\tau(L)\cap \pi_\alpha=R_2(q)$ in the latter. 
\end{proof}

\begin{lemma}\label{l:comp} The following hold.
\begin{enumerate}
 \item $(\overline \Phi)_{2'}$ centralizes $(\Outdiag L)_2$.
 \item If $L\neq O_8^+(q)$, then $(\Outdiag L)_{2'}\rtimes (\overline \Phi)_{2'}$ is a normal $2$-complement in $\Out L$.
 \item Let $L\neq O_8^+(q)$ and let $S$ be a Sylow $2$-subgroup of $\overline G$. If two elements of $S$ are conjugated in $\overline G$, then they are conjugated in $S$.
\end{enumerate}

\end{lemma}

\begin{proof} 
It is clear that (iii) follows from (ii), while (ii) follows from (i) because $\overline \Phi\overline \Gamma$ normalizes $(\Outdiag L)_{2'}\rtimes (\overline \Phi)_{2'}$. So it sufficient to prove (i). 

We may assume that $|(\Outdiag L)_2|>2$. If $\Outdiag L$ is not cyclic, then $\overline \Phi$ centralizes $\Outdiag L$ (see \cite[Theorem 2.5.12(h)]{98GorLySol}). Thus we are left with  $L_n^\varepsilon(q)$, where $(n,q-\varepsilon)_2>2$, and $O_{2n}^\varepsilon(q)$, where $n$ is odd and $(q-\varepsilon, 4)=4$. Let $q_1=p^{(f)_{2}}$. We can choose generators $\eta$ and $\phi$ of $(\Outdiag L)_2$ and $(\Phi)_{2'}$, respectively, so that $\eta^\phi=\eta^{q_1}$ if $\varepsilon=+$ and $\eta^\phi=\eta^{q_1^2}$ if $\varepsilon=-$. Since $|\eta|$ divides $(q-\varepsilon)_2$ and $(q-\varepsilon)_2=(p^f-\varepsilon)_2=(p^{(f)_{2}}-\varepsilon)_2=(q_1-\varepsilon)_2,$ the claim follows.
\end{proof}

We are ready to describe all almost simple $W_r$-groups (recall Notation \ref{n:auto}).

\begin{theorem} \label{t:1} Let $r$ be an odd prime, let $L$ be a simple group of Lie type, $L\leq G\leq \Aut L$ and $\overline G=G/L$.  Then $G$ is $W_r$-group if and only if one of the following holds:
\begin{enumerate}
\item $L$ is a $W_r$-group and $2, r\not\in\pi(\overline G)$;
 \item $L\neq O_8^+(q)$, $2\in\pi(\overline G)$ and $\overline G$ is conjugate to a subgroup of $$((\Outdiag L)_{2'}\rtimes (\overline \Phi)_{\{2,r\}'})\rtimes \langle \theta\rangle=((\Outdiag L)_{2'}\rtimes \langle \theta\rangle) \rtimes (\overline \Phi)_{\{2,r\}'},$$ where $1\neq \theta \in  (\Outdiag L)_{2}\rtimes (\overline \Phi\overline \Gamma)_{2}$ and the only involution $\alpha$ of $\langle \theta\rangle$ is as in Tables {\rm \ref{tab:2even}} and {\rm \ref{tab:2odd}}.
 
 \item $L=O_8^+(q)$, $q$ is a square, $r\in R_6(q)$, $2\in\pi(\overline G)$ and $\overline G$ is conjugate to a subgroup of $$(\langle \tau \rangle\times (\overline \Phi)_{\{2,r\}'})\rtimes \langle \theta\rangle=(\langle \tau \rangle \rtimes \langle \theta\rangle)\times (\overline \Phi)_{\{2,r\}'},$$ where $\theta \in (\langle \gamma\rangle\times (\overline \Phi)_{2})\setminus \langle\gamma\rangle$, or $q$ is odd and $\overline G$ is conjugate to a subgroup of $$(\Outdiag L\rtimes \langle \gamma\rangle)\times  \overline \Phi_{r'} $$
 trivially intersecting $\Outdiag L\rtimes  \langle \gamma\rangle$.
\end{enumerate}


\end{theorem}

\begin{proof}
By Lemma \ref{l:cr}, if $|\overline G|$ is odd, then $G$ is a $W_r$-group if and only if (i) holds. So we may assume that $|\overline G|$ is even. If (ii) or (iii) holds, then $G$ is a~$W_r$-group by Lemma \ref{l:cr}(ii) and Propositions \ref{p:2even} and \ref{p:2odd}.

Conversely, suppose that $G$ is $W_r$-group. By Lemma \ref{l:general}, Sylow $2$-subgroups of $\overline G$ are cyclic. By \cite[Corollary 5.14]{08Isa}, it follows that $\overline G$ has a normal 2-complement.

If $r=3$, then (ii) holds by Lemma \ref{l:w3}. Indeed, the lemma states that $O^{\{2,3\}'}(G)\leq PGL_2(3^f)$, or $f$ is even and $O^{\{2,3\}'}(G)\leq PGL_2^*(3^f)$, or equivalently, that  $\overline G\leq \langle\alpha\rangle \times (\overline\Phi)_{\{2,3\}'}$, where $\alpha=\delta$ or $\delta\varphi$. Such $\alpha$ is as in Table \ref{tab:2odd}. 


Now let $r>3$ and assume first that $L\neq O_8^+(q)$. Applying Lemma \ref{l:comp}(ii) and noting that  $r$ does not divide $|\Outdiag L|$ by Lemma \ref{l:general}(i), we see that $(\Outdiag L)_{2'}\rtimes (\overline \Phi)_{\{2,r\}'}$ is the normal $\{2,r\}'$-Hall subgroup of $\Out L$.  Thus $\overline G$ is conjugate to a subgroup of $((\Outdiag L)_{2'}\rtimes (\overline \Phi)_{\{2,r\}'})\rtimes \langle \theta\rangle$, where $1\neq \theta \in (\Outdiag L)_{2}\rtimes (\overline \Phi\overline \Gamma)_{2}$. By Propositions \ref{p:2even} and \ref{p:2odd}, the only involution of $\langle \theta\rangle$ is conjugate to $\alpha$, where $\alpha$ is as in Table \ref{tab:2even} or Table \ref{tab:2odd}. In view of Lemma \ref{l:comp}(iii), we may assume that this involution is equal to $\alpha$. It remains to observe that $\theta$ centralizes  $(\overline \Phi)_{\{2,r\}'}$   
by Lemma \ref{l:comp}(i), and so the equality in (ii) holds.

Let $L=O_8^+(q)$.  It is clear that $\overline G\leq S\times \overline \Phi_{r'}$, where $S=\Outdiag L\rtimes \overline \Gamma$. Note that $S\simeq S_3$ if $q$ is even and $S\simeq S_4$ if $q$ is odd.
Let $\theta$ be a generator of a Sylow $2$-subgroup of $\overline G$. 
By Tables \ref{tab:2even} and \ref{tab:2odd}, it follows that $f$ is even, $r\in R_6(q)$ and $\langle \theta\rangle \cap S=1$. 

Suppose that the projection $\overline G_1$ of the Sylow $3$-subgroup of $\overline G$ on $S$ is not trivial. Since $\overline G$ has a normal 2-complement, the projection $\theta_1$ of $\theta$ on $S$ must normalizes $\overline G_1$. Since  only transpositions and the identity normalize a group of order $3$ in both $S_3$ and $S_4$, it follows that, up to conjugacy,  $\overline G_1=\langle \tau\rangle$ and $\theta_1=\gamma$ or $1$. Hence $\overline G\leq(\langle \tau \rangle \rtimes \langle \theta\rangle)\times (\overline \Phi)_{\{2,r\}'}$, where $\theta \in (\langle \gamma\rangle\times (\overline \Phi)_{2})\setminus \langle\gamma\rangle$. 

If $\overline G_1=1$, then, up to conjugacy, $\overline G\leq (\Outdiag L\rtimes \langle \gamma\rangle)\times  \overline \Phi_{r'}$. If $q$ is even, then $\overline G$ is as in the previous case, so here we may assume that $q$ is odd. The condition  $\langle \theta\rangle \cap S=1$ is equivalent to $\overline G\cap (\Outdiag L\rtimes  \langle \gamma\rangle)=1$.
\end{proof}

We conclude this section with two corollaries of Theorem \ref{t:1}, demonstrating how to apply this theorem to a particular $L$ or a particular $r$. The first corollary describes almost simple $W_r$-groups with socle $L_2(q)$. 


\begin{cor}\label{c:l2}
Let $r$ be an odd prime, $L=L_2(q)$,  where $q=p^f$, and $L\leq G\leq \Aut L$. Set $\overline G=G/L$.
Then $G$ is a $W_r$-group if and only if one of the following holds:
\begin{enumerate}
\item $r=p$, end either $\overline G\leq \langle \delta\rangle \times (\overline \Phi)_{\{2,r\}'}$, or $f$ is even and $\overline G\leq \langle \delta\varphi\rangle \times (\overline \Phi)_{\{2,r\}'}$;
\item $r\in \pi(q-1)$, $q\not\equiv 1\pmod 4$, and $\overline G\leq  (\overline \Phi)_{\{2,r\}'}$;
 \item $r\in \pi(q+1)$, $q$ is even, and $\overline G\leq (\overline \Phi)_{r'}$;
 \item $r\in \pi(q+1)$, $q\equiv 1\pmod 4$, $\overline G\leq \langle \delta\rangle \times (\overline \Phi)_{r'}$ and $\delta\not\in \overline G$.
\end{enumerate}
\end{cor}

\begin{proof}
By \cite[Tables 4,6]{05VasVd.t}, it follows that $L$ is a $W_r$-group if and only if either $r=p$, or  
$r$ divides $q-\epsilon$ for some $\epsilon\in\{+1,-1\}$ and $(q-\epsilon)_2\leq 2$.

Observe that $\Out L$ is abelian. By Theorem \ref{t:1}, the group $G$ is a $W_r$-group if and only if $\overline G\leq \langle \theta\rangle \times (\overline\Phi)_{\{2,r\}'}$, where $\theta\in \Outdiag L\times (\overline \Phi)_2$, and the only involution $\alpha$ of $\langle \theta\rangle$ is as in Table \ref{tab:2even} or Table \ref{tab:2odd}. We may assume that $\alpha\in\overline G$.  In particular, we may assume that $r\not\in R_1(q)$ since $R_1(q)$ does not appear in  Tables \ref{tab:2even} and \ref{tab:2odd}. 

Let $r=p$. By Table \ref{tab:2odd}, we have $\alpha=\delta$, or $f$ is even and $\alpha=\delta\varphi$. In either case, $\theta=\alpha$ and so $\overline G$ is as stated in (i). 

Let $r\in R_2(q)$. By Tables \ref{tab:2even} and \ref{tab:2odd},  either $q$ is even, in which case $\alpha=\varphi$ and so $\overline G\leq (\overline \Phi)_{r'}$, or $q$ is an odd square and $\alpha\in\{\varphi, \delta\varphi\}$, in which case $\overline G$ is as stated in (iv).  
\end{proof}
  
The second corollary describes almost simple $W_5$-groups of Lie type. After the proof, we check this description against the description of almost simple $W_5$-groups given in \cite[Theorem 1]{24GGKN.t}. Since \cite[Theorem 1]{24GGKN.t} imposes conditions on $O^{\{2,5\}'}(G)$ rather then on $G$, we also impose conditions on $O^{\{2,5\}'}(G)$ (it is clear that $O^{\{2,5\}'}(\overline G)=\overline{O^{\{2,5\}'}(G)}$).

\begin{cor}\label{c:w5}
Suppose that $L$ is a simple group of Lie type and  $L\leq G\leq \Aut L$. Then $G$ is a $W_5$-group if and only if one of the following holds:
\begin{enumerate}
\item $L=L_2(5^f)$ and either $O^{\{2,5\}'}(\overline G)\leq \langle \delta\rangle$,  or $f$ is even and $O^{\{2,5\}'}(\overline G)=\langle \delta\varphi\rangle$;
 \item $L=L_2(q)$, where $q\equiv \epsilon\pmod 5$ for some $\epsilon\in\{+1,-1\}$, $(q-\epsilon)_2\leq 2$, and either $O^{\{2,5\}'}( \overline G)=1$, or $\epsilon=-1$, $q$ is a square, and
 $O^{\{2,5\}'}(\overline G)=\langle \varphi\rangle$, or $q$ is an odd square and $O^{\{2,5\}'}(\overline G)=\langle \delta\varphi\rangle$; 
\item $L=L_3^\varepsilon(q)$, where $q\equiv -\varepsilon\pmod 5$, $q$ is even, and  either $O^{\{2,5\}'}(\overline G)=1$, or $\varepsilon=+$, $q$ is a~square,  $O^{\{2,5\}'}(\overline G)$ is a subgroup of even order in $\Outdiag L\rtimes \langle \alpha\rangle$ for $\alpha\in\{\varphi, \varphi\gamma\}$;
\item $L={}^2B_2(q)$, where $q>2$, and $O^{5'}( \overline G)=1$;
 \item $q\equiv \pm 2\pmod 4$, $O^{\{2,5\}'}( \overline G)=1$, and either $L=L_4^\varepsilon(q)$ with $(q+\varepsilon)_2\leq 2$, or $L=L_5^\varepsilon(q)$ with $q$ even, or $L=S_4(q)$ with $q>2$. 
\end{enumerate}
\end{cor}

\begin{proof}
 We may assume that $5\in\pi(L)$. Let $q=p^f$. Note that $5\in\{p\}\cup R_1(q)\cup R_2(q)\cup R_4(q)$. Furthermore, if $5\in R_4(q)$, then $f$ is odd and if $5\in R_2(q)$, then $(f)_2\leq 2$. 
 
If $L=L_2(q)$, then the assertion follows from Corollary \ref{c:l2}. If $L\neq L_2(q)$, then we may assume that $p\neq 5$ by \cite[Tables 4--7]{05VasVd.t}.

Let $L=L_3^\varepsilon(q)$. Then $5\in R_1(q)\cup R_2(q)$, and by \cite[Tables 4,6]{05VasVd.t}, we have $10\not\in\omega(L)$ if and only if $q$ is even and $5\in R_2(\varepsilon q)$. By Table \ref{tab:2even}, if $G$ is as stated in (iii), then $G$ is a $W_5$-group. Conversely, suppose that  $G$ is a $W_5$-group. By Theorem \ref{t:1}, if $O^{\{2,5\}'}( \overline G)\neq 1$, then $O^{\{2,5\}'}(\overline G)$ is conjugate to a~subgroup of $H=\Outdiag L\rtimes \langle \theta\rangle$, where $1\neq \theta \in  (\overline{\Phi\Gamma})_2$ and the only involution $\alpha$ of $\langle \theta\rangle$ is as in Table \ref{tab:2even}.  Since $H$ is normal in $\Out L$, we have $\overline G\leq H$. By Table \rm \ref{tab:2even} we see that $\varepsilon=+$, $f$ is even and $\alpha\in\{\varphi, \varphi\gamma\}$. Furthermore, $(f)_2=2$, and hence $\theta=\alpha$.

If $L={}^2B_2(q)$, then $L$ is a $W_5$-group for all $q$ and $|\Out L|$ is odd. So $G$ is a $W_5$-group if and only if $O^{5'}( \overline G)=1$.

If $L\neq L_2(q)$, $L_3^\varepsilon(q)$, $^2B_2(q)$, then it follows from \cite[Tables 4--7]{05VasVd.t} that divisors of $q^2-1$ are adjacent to $2$ in $GK(G)$ (see also Table \ref{tab:ex} and Lemma \ref{l:I} in Section \ref{s:t2lie}), and therefore $5\in R_4(q)$, or equivalently, $q\equiv \pm 2\pmod 5$. Furthermore, elements of $R_4(q)$ are not adjacent to $2$ in $GK(L)$ if and only if $L$ is as specified in (v). Since $f$ is odd, the set $R_4(q)$ does not appear in Tables \ref{tab:2even} and \ref{tab:2odd}. Thus $G$ is a $W_5$-group if and only if $2,5\not\in\pi(\overline G)$.
\end{proof}

\begin{remark}\label{rem:w5} It is easily seen that (i), (ii) and (iv) of Corollary \ref{c:w5} are equivalent to (1),  (2) and (7) of \cite[Theorem 1]{24GGKN.t}, respectively, while (v) of Corollary \ref{c:w5} incorporates (4)-(6) of \cite[Theorem 1]{24GGKN.t}.  Also (iii) of Corollary \ref{c:w5} is equivalent to (3) of \cite[Theorem 1]{24GGKN.t} but this requires some explanations in the case when $O^{\{2,5\}'}(\overline G)$  contains a conjugate of $\varphi\gamma$.

Suppose that $L=L_3^\varepsilon(q)$, where $q\equiv -\varepsilon\pmod 5$, $q$ is even and $O^{\{2,5\}'}(\overline G)$ contains a~conjugate of $\varphi\gamma$. Then $\varepsilon=+$ and $q$ is a square, and
(iii) states that $O^{\{2,5\}'}(\overline G)\leq \Outdiag L\rtimes \langle \varphi\gamma\rangle$, while (3) of \cite[Theorem 1]{24GGKN.t} states that $O^{\{2,5\}'}(\overline G)=\langle \varphi\gamma\rangle$. As we noticed, $(f)_2=2$ and then an easy calculation shows that $\varphi\gamma$ centralizes $\Outdiag L$. Thus 
$O^{\{2,5\}'}(\overline G)\leq \Outdiag L\rtimes \langle \varphi\gamma\rangle$ implies $O^{\{2,5\}'}(\overline G)=\langle \varphi\gamma\rangle$.
\end{remark}

\section{2-Independence numbers for almost simple groups of Lie type}\label{s:t2lie}

We keep notation of the previous section. Recall that $L={}^d\Sigma(q)$, where $q=p^f$, 
and $L\leq G\leq \Aut L$. Our aim is to find $\tau(G)$ and $t(2,G)$ for all such $G$.
Suppose that $\rho$ is a $2$-coclique in $GK(G)$ and $\rho'=\rho\setminus\{2\}$. It follows from Lemma~\ref{l:simple} that $\rho'\cap \pi(\overline G)=\varnothing$. In particular, $\rho'\subseteq \tau(L)$ and $t(2, G)\leq t(2,L)$. Also recall that by Lemma \ref{l:partition}, we have a~partition of $\tau(L)$ into cocliques and we will write this partition in the form $\sqcup \tau(L)_i$.

\begin{lemma}\label{l:general_ind} 
Let $\tau(L)=\sqcup \tau(L)_i$ be the partition of $\tau(L)$ with respect to $\sim$. 
Then $t(2,G)=|I(G)|+1$, where $I(G)=\{i\mid \tau(G)\cap \tau(L)_i\neq \varnothing\}$. Furthermore, the following hold.

\begin{enumerate}
  \item If $|\overline G|$ is odd,  then $\tau(G)\cap \tau(L)_i=\tau(L)_i\setminus \pi(\overline G)$. 
  \item If Sylow $2$-subgroups of $\overline G$ are cyclic and $\alpha$ is some involution of $\overline G$, then $\tau(G)\cap \tau(L)_i=\tau(L)_i\setminus(\pi(\overline G)\cup \pi_\alpha))$. 
\end{enumerate}
\end{lemma}

\begin{proof} This follows from Lemmas \ref{l:partition} and \ref{l:cr}.
\end{proof}

Thus our task is to determine when the set $\tau(L)_i\setminus \pi(\overline G)$ or $\tau(L)_i\setminus(\pi(\overline G)\cup \pi_\alpha)$ is not empty. If $t(2,G)=t(2,L)$ and $|\overline G|$ is odd, then the partition of $\tau(G)$ is clear. If $t(2,G)<t(2,L)$ or $|\overline G|$ is even, then we will provide the non-empty sets $\tau(G)_i=\tau(L)_i\cap \tau(G)$. 
For brevity, if $G$ is specified and $\pi\subseteq\pi(G)$, then we write $$\pi^*=\pi\setminus \pi(\overline G).$$

 Expecting \cite[Tables 4--7]{05VasVd.t}, one can observe that, in general, the elements of $\tau(L)/{\sim}$ are sets of the form  $R_i(q)$, but there are some exceptions. These exceptions are given in Table~\ref{tab:ex}, and for the other groups, the properties of $\tau(L)/{\sim}$ are written in Lemma \ref{l:I}.

\begin{table}[htb]
\caption{Simple groups $L$ of Lie type with exceptional $\tau(L)$}\label{tab:ex}
\begin{center}
\begin{tabular}{|l|l|l|}
\hline
$L$& Conditions on $L$& $\tau(L)/{\sim}$\\
\hline
$L_2(q)$& $5\leq q\equiv \varepsilon\pmod 4$ & $\{p\}$, $R_2(\varepsilon q)$ \\
$L_3(q)$&$q\geq 4$ even, $(q-1)_3=3$& $\{3\}$, $R_2(q)$, $R_3(q)$\\
$U_3(q)$&$q\geq 4$ even, $(q+1)_3=3$& $\{3\}$, $R_1(q)$, $R_6(q)$\\
${}^2B_2(q)$&$q=2^f$, $f\geq 3$ odd& $R_1(q)$, $\pi(q\pm \sqrt{2q}+1)$\\ 
${}^2G_2(q)$&$q=3^f$, $f\geq 3$ odd& $\pi(q\pm\sqrt{3q}+1)$\\ 
${}^2F_2(q)$&$q=2^f$, $f\geq 3$ odd& $R_6(q)$, $\pi(q^2\pm\sqrt{2q^3}+q\pm\sqrt{2q}+1)$\\ 
${}^2F_4(2)'$&&$\{13\}$\\
\hline
\end{tabular}
\end{center}
\end{table}

\begin{lemma} \label{l:I} Suppose that $L$ is not a group in Table {\rm \ref{tab:ex}}.
Then we can choose $I\subseteq \mathbb N$ such that $\tau(L)/{\sim}\,=\{R_i(\epsilon q)\mid i\in I\}$, where $\epsilon=\varepsilon$ if $L=L_n^\varepsilon(q)$ with $n\geq 3$,  $O_{2n}^\varepsilon(q)$ with $n$ odd, or $E_6^\varepsilon(q)$, and $\epsilon=+$ otherwise. Furthermore, if $I\cap\{1,2\}\neq \varnothing$, then $L=L_2(q)$ or $L_3^\varepsilon(q)$.
\end{lemma}

In the next three lemmas, we find $\tau(G)$ and $t(2,G)$ in the case when $L$ has exceptional $\tau(L)$. Note that the group $^2F_4(2)'$ is not a Ree group. 

\begin{lemma}\label{l:ind_sr}
Let $L$ be a Suzuki or Ree group. Then either $t(2,G)=t(2,L)$, or $L={}^2B_2(32)$, $G=\Aut L$,  $t(2,G)=t(2,L)-1=3$ and $\tau(G)=\{31\}\sqcup\{41\}$.
\end{lemma}

\begin{proof}
 This is \cite[Theorem 3]{23Gre.t} but we repeat the proof here for completeness. Since $\pi(\overline G)\subseteq \pi(f)$ and $f$ is odd, we apply Lemma \ref{l:general_ind}(i). Taking $\tau(L)_i$ from Table \ref{tab:ex} and using Lemmas \ref{l:bz0} and \ref{l:bzsr0}, we see that 
 $\tau(L)_i\setminus\pi(f)\neq \varnothing$ for all $i$, and so $\tau(L)_i\setminus\pi(\overline G)\neq\varnothing$ for all $i$, provided that $L\neq{}^2B_2(32)$, then  If $L={}^2B_2(32)$ and $G=\Aut L$,  then $\tau(L)=\{5\}\sqcup\{31\}\sqcup\{41\}$ and the claim follows.
\end{proof}

\begin{lemma}\label{l:ind_l2}
Let $L=L_2(q)$, where $q=p^f\geq 5$ is odd and $q\equiv \varepsilon\pmod 4$. Suppose that Sylow $2$-subgroups of $\overline G$ are cyclic.
\begin{enumerate}
\item If $|\overline G|$ is odd or $\delta\varphi\in\overline G$, then $t(2,G)=2$ or $3$ depending on whether $p$ divides $|\overline G|$ or not,
with $\tau(G)=R_2(\varepsilon q)^*$ or $\tau(G)=\{p\}\sqcup R_2(\varepsilon q)^*$, respectively.

\item If $\varphi \in \overline G$, then $t(2,G)=2$ and $\tau(G)=R_2(\varepsilon q)^*$.

\item If $\delta\in \overline G$, then $t(2,G)=1$ or $2$ depending on whether $p$ divides $|\overline G|$ or not, with $\tau(G)=\{p\}$ in the latter case.
\end{enumerate}

\end{lemma}

\begin{proof}
Note that $R_2(\varepsilon q)\setminus\pi(f)\neq \varnothing$ by Lemma \ref{l:bz0}. In particular, the claim is immediate if $|\overline G|$ is odd. If $|\overline G|$ is even, then the claim follows from Table \ref{tab:2odd}.
\end{proof}

\begin{lemma}\label{l:ind_l3}
Let  $L=L^\varepsilon_3(q)$, where $q\geq 4$ is even and $(q-\varepsilon)_3=3$, and $L\leq G\leq\Aut L$.  Suppose that Sylow $2$-subgroups of $\overline G=G/L$ are cyclic. 

\begin{enumerate}
 \item If $2,3\not\in\pi(\overline G)$, then $t(2,G)=t(2,L)=4$.
 \item If $2\not\in\pi(\overline G)$ and $3\in\pi(\overline G)$, then $t(2,G)=3$, and $\tau(G)=R_2(\varepsilon q)^*\sqcup R_3(\varepsilon q)^*$.
 \item If $\overline G$ contains a conjugate of $\gamma$, then $t(2,G)=2$ and $\tau(G)=R_3(\varepsilon q)^*$.
 \item If $\varepsilon=+$ and $\overline G$ contains a conjugate of $\varphi\gamma$, then $t(2,G)=3$ and $\tau(G)=R_2(q)^*\sqcup R_3(q_0)^*$.
\item If $\varepsilon=+$ and $\overline G$ contains a conjugate of $\varphi$, 
then $t(2,G)=3$ and $\tau(G)=R_2(q)^*\sqcup R_6(q_0)^*$, unless $L=L_3(4)$, in which case $t(2,G)=2$ and $\tau(G)=\{5\}$.
\end{enumerate}

\end{lemma}

\begin{proof} Observe that $3\in R_1(\varepsilon q)$, so $R_i(\varepsilon q)^*=R_i(\varepsilon q)\setminus\pi(\overline G)=R_i(\varepsilon q)\setminus\pi(f)$ for $i=2,3$. Since $q\neq 8$,  it follows from Lemma \ref{l:bz0}, that both $R_2(\varepsilon q)^*$ and $R_3(\varepsilon q)^*$ are not empty. So (i) and (ii) are immediate. 

Let $2\in\pi(\overline G)$. If $\overline G$ contains a conjugate of $\gamma$, then it follows from Table \ref{tab:2even} that the only non-empty part of $\tau(G)$ is contained in $R_3(\varepsilon q)$, and we have (iii). So we may assume that $\varepsilon=+$ and $f$ is even. Using Table \ref{tab:2even} and Lemma \ref{l:bz0}(ii), we complete the proof.
\end{proof}

We are ready to prove the main result of this section.

\begin{theorem}\label{t:2} Let $L$ be a simple group of Lie type, $L\leq G\leq \Aut L$ and set $\overline G=G/L$.
\begin{enumerate}
 \item If Sylow $2$-subgroups of $\overline G$ are not cyclic then $t(2,G)=1$.
 \item If $|\overline G|$ is odd, then either $t(2,G)=t(2,L)$ and $\tau(G)=\sqcup \tau(L)_i^*$, or one of the following holds:
 \begin{enumerate}
 \item $PGL_3^\varepsilon(q)\leq G$, where $q\geq 4$ is even and $(q-\varepsilon)_3=3$,     $t(2,G)=3$ and $\tau(G)=R_2(\varepsilon q)^*\sqcup R_3(\varepsilon q)^*$; 
  \item $G=\Aut (L_2(8))$, $t(2,G)=2$,   $\tau(G)=\{7\}$;
  \item $G=L_3(8).3$, $t(2,G)=2$, $\tau(G)=\{73\}$;
  \item $G=\Aut({}^2B_2(32))$, $t(2,G)=3$ and $\tau(G)=\{31\}\sqcup\{41\}$. 
 
\end{enumerate}
  \item If Sylow $2$-subgroups of $\overline G$ are non-trivial and cyclic, then $t(2,G)$ and $\tau(G)$ are as in Tables {\rm \ref{tab:adj_even}, \ref{tab:adj_odd}} and {\rm \ref{tab:adj_ex}}.  
\end{enumerate}
\end{theorem}

\begin{remark}
In Tables  {\rm \ref{tab:adj_even}, \ref{tab:adj_odd}} and {\rm \ref{tab:adj_ex}}, we write ``$\alpha\in \overline G$'' in the third column to indicate that $\overline G$ contains some conjugate of $\alpha$. Similarly, ``$\alpha\not\in \overline G$'' indicates that ``$G$ does not contain conjugates of $\alpha$''.
\end{remark}

\begin{table}[t]
\caption{2-Independence numbers for finite almost simple classical groups $G$  in characteristic $2$ with $G/\Soc(G)$ having a non-trivial cyclic Sylow $2$-subgroup}\label{tab:adj_even}

\renewcommand{\arraystretch}{1.3}
\begin{tabular}{|l|l|l|l|l|}
\hline
$L$&Conditions on $L$& Conditions on $\overline G$& $t(2,G)$& $\tau(G)$\\
\hline
$L_n(q)$&$n=2$&&2&$R_2(q)^*$\\
\cline{2-5}
&$n=3$, $q=4$& $\varphi\in\overline G$&2&$\{5\}$\\
&&$\varphi\gamma\in\overline G$&3&$\{5\}$,$\{7\}$\\
&&$\gamma\in \overline G$&2&$\{7\}$\\
\cline{2-5}
&$n\geq 3$ odd, $(n,q)\neq (3,4)$& $\varphi\in\overline G$&3&$R_{n-1}(q)^*$, $R_{2n}(q_0)^*$\\
&&$\varphi\gamma\in\overline G$&3&$R_{n-1}(q)^*$, $R_{n}(q_0)^*$\\
&&$\gamma\in \overline G$&2&$R_n(q)^*$\\
\cline{2-5}
&$n=4$, $q=4$& $\varphi\in\overline G$&2&$\{17\}$\\
&&$\varphi\gamma\in\overline G$&3&$\{7\}$,$\{17\}$\\
&&$\gamma\in \overline G$&2&$\{7\}$\\
\cline{2-5}
&$n\geq 4$ even, $(n,q)\neq (4,4)$& $\varphi\in\overline G$&3&$R_n(q)^*$, $R_{2n-2}(q_0)^*$\\
&& $\varphi\gamma\in\overline G$&3&$R_n(q)^*$, $R_{n-1}(q_0)^*$\\
&&$\gamma\in \overline G$&2&$R_{n-1}(q)^*$\\
\hline
$U_n(q)$&$n\geq 3$ odd&&2& $R_{n}(-q)^*$\\
\cline{2-5}
&$n\geq 4$ even&&2& $R_{n-1}(-q)^*$\\
\hline
$S_{2n}(q)$&$n\geq 3$ or $f$ even&&2&$R_{2n}(q)^*$\\
\cline{2-5}
&$n=2$, $f$ odd&&1&\\
\hline
$O_{2n}^+(q)$
&$n\geq 4$ even&$\gamma\not\in\overline G$&2&$R_{2n-2}(q)^*$\\
&&$\gamma\in \overline G$&1&\\
\cline{2-5}
&$n\geq 5$ odd&$\varphi\in\overline G$&3&$R_{2n-2}(q)^*$, $R_{2n}(q_0)^*$\\
&& $\varphi\gamma\in\overline G$&3&$R_{2n-2}(q)^*$, $R_{n}(q_0)^*$\\
&&$\gamma\in \overline G$&2&$R_n(q)^*$\\
\hline 
$O_{2n}^-(q)$&&&2& $R_{2n}(q)^*$\\
\hline
\end{tabular}
\end{table}

\begin{table}
\caption{2-Independence numbers for finite almost simple classical groups $G$  in odd characteristic with $G/\Soc(G)$ having a non-trivial cyclic Sylow $2$-subgroup}\label{tab:adj_odd}

\renewcommand{\arraystretch}{1.3}
\begin{tabular}{|l|l@{\hspace{0pt}}|l|l|l|}
\hline
$L$&Conditions on $L$& Conditions on $\overline G$& $t(2,G)$& $\tau(G)$\\
\hline
$L_n(q)$&$n=2$&$\varphi\in \overline G$&2&$R_2(q)^*$\\
&&$\varphi\delta\in\overline G$, $p\nmid |\overline G|$&3&$\{p\}$, $R_2(q)^*$\\
&&$\varphi\delta\in\overline G$, $p\mid |\overline G|$&2&$R_2(q)^*$\\
&&$\delta\in \overline G$, $p\nmid |\overline G|$&2&$\{p\}$\\
&&$\delta\in \overline G$, $p\mid |\overline G|$&1&$\{p\}$\\
\cline{2-5}
&$n\geq 3$ odd&$\varphi\in\overline G$&2&$R_{2n}(q_0)^*$\\
&&$\varphi\gamma\in\overline G$&2&$R_{n}(q_0)^*$\\
&&$\gamma\in \overline G$&2&$R_{n}(q)^*$\\
\cline{2-5}
&\multirow{10}{*}{\begin{tabular}{@{\hspace{0pt}}l|l}$n\geq 4$ even&$(n)_2>(q-1)_2$\\
&\\
&\\
&\\
\cline{2-2}
&$(n)_2=(q-1)_2$\\
&\\
&\\
&\\
\cline{2-2}
&$(n)_2<(q-1)_2$\\
&\\
\end{tabular}}&
$\varphi\in\overline G$&2&$R_{2n-2}(q_0)^*$\\
&& $\varphi\gamma\in\overline G$&2& $R_{n-1}(q_0)^*$\\
&&$\gamma\in \overline G$ or $\gamma'\in \overline G$&2&$R_{n-1}(q)^*$\\
&&$\delta\in \overline G$&1&\\
\cline{3-5}
&&$\varphi\in\overline G$&3&$R_n(q)^*$, $R_{2n-2}(q_0)^*$\\
&& $\varphi\gamma\in\overline G$&3&$R_n(q)^*$, $R_{n-1}(q_0)^*$\\
&&$\gamma\in \overline G$ or $\gamma'\in \overline G$&2&$R_{n-1}(q)^*$\\
&&$\delta\in \overline G$&1&\\
\cline{3-5}
&&
$|(\Outdiag L\rtimes \overline \Gamma)\cap \overline G|$ odd &2&$R_n(q)^*$\\
&&$|(\Outdiag L\rtimes \overline \Gamma)\cap \overline G|$ even &1&\\
\hline
$U_n(q)$&$n\geq 3$ odd&&2& $R_n(-q)^*$\\
\cline{2-5}
&$(n)_2\geq (q+1)_2$&$\gamma\in \overline G$ or $\gamma'\in \overline G$
&2&$R_{n-1}(-q)^*$\\
\cline{2-5}
&$2=(n)_2<(q+1)_2$&$\gamma'\in\overline G$&2& $R_{n}(-q)^*$\\
\cline{2-5}
&\multicolumn{2}{c|}{otherwise}&1&\\
\hline
$S_{2n}(q)$,&&$\delta\not\in \overline G$&2&$R_{2n}(q)^*$\\
\cline{2-5}
$O_{2n+1}(q)$&&$\delta\in \overline G$&1&\\
\hline
$O_{2n}^+(q)$
&$n\geq 4$ even&$|(\Outdiag L\rtimes \overline \Gamma)\cap \overline G| $ odd&2&$R_{2n-2}(q)^*$\\
\cline{2-5}
&$n\geq 5$ odd, $q\not\equiv 1\pmod 8$&$\gamma\in \overline G$ or $\gamma'\in\overline G$&2&$R_n(q)^*$\\
\cline{2-5}
&\multicolumn{2}{c|}{otherwise}&1&\\
\hline 
$O_{2n}^-(q)$&$n\geq 4$ even or $q\not\equiv 7\pmod 8$&$\delta\not\in\overline G$&2&$R_{2n}(q)^*$\\
\cline{2-5}
&\multicolumn{2}{c|}{otherwise}&1&\\
\hline
\end{tabular}
\end{table}

\begin{table}
\caption{2-Independence numbers for finite almost simple exceptional groups $G$ with $\overline G$ having a non-trivial cyclic Sylow $2$-subgroup}\label{tab:adj_ex}

\renewcommand{\arraystretch}{1.3}
\begin{tabular}{|l|l|l|l|l|}
\hline
$L$&Conditions on $L$& Conditions on $\overline G$& $t(2,G)$& $\tau(G)$\\
\hline
$G_2(q)$, $q\geq 3$&$p=3$, $f$ odd& &2&$R_{3}(q)^*$\\
\cline{2-5}
&$p\neq 3$ or $f$ even&&2&$R_{6}(q)^*$\\
\hline
$F_4(q)$&$q$ even, $f$ even&&3&$R_{8}(q)^*$, $R_{12}(q)^*$\\
\cline{2-5}
&$q$ even, $f$ odd&&2&$R_{8}(q)^*$\\
\cline{2-5}
&$q$ odd&&2&$R_{12}(q)^*$\\
\hline
$^2F_4(2)'$&&&$2$&$\{13\}$\\
\hline 
$^3D_4(q)$&&&2& $R_{12}(q)^*$\\
\hline
$E_6(q)$&$q$ even &$\varphi\in \overline G$&4& $R_{8}(q)^*$, $R_{12}(q)^*$, $R_{18}(q_0)^*$\\
&&$\varphi\gamma\in \overline G$&4& $R_{8}(q)^*$, $R_{12}(q)^*$, $R_{9}(q_0)^*$\\
& &$\gamma\in \overline G$&2& $R_{9}(q)^*$\\
\cline{2-5}
&$q$ odd&$\varphi\in \overline G$&3&  $R_{12}(q^*$, $R_{18}(q_0)^*$\\
&& $\varphi\gamma\in \overline G$&3&  $R_{12}(q)^*$, $R_{9}(q_0)^*$\\
& &$\gamma\in \overline G$&2& $R_{9}(q)^*$\\
\hline
$^2E_6(q)$& & &2& $R_{18}(q)^*$\\
\hline
$E_7(q)$&$q$ even & &3& $R_{14}(q)^*$, $R_{18}(q)^*$\\
\cline{2-5}
&$q\equiv 1\pmod 4$ &$\delta\not\in \overline G$ &3& $R_{14}(q)^*$, $R_{18}(q)^*$\\& &$\delta\in \overline G$ &1& \\
\cline{2-5}
&$q\equiv 3\pmod 4$ & &1& \\
\hline 
$E_8(q)$& &  & $4$&  $R_{20}(q)^*$, $R_{24}(q)^*$, $R_{30}(q)^*$\\
\hline
\end{tabular}
\end{table}

\begin{proof}
By Lemma \ref{l:general}(ii), we have (i). 

If $L={}^2F_4(2)'$, then $\tau(G)=\tau(L)=\{13\}$ by \cite{85Atlas}, as stated in Table \ref{tab:adj_ex}. If $L$ is a Suzuki--Ree group, or $L_2(q)$, where $q$ is odd, or $L_3^\varepsilon(q)$, where $q$ is even and $(q-\varepsilon)_3=3$,
then the claim follows from Lemmas \ref{l:ind_sr}, \ref{l:ind_l2} and \ref{l:ind_l3}, respectively. In the case $L=L_3^\varepsilon(q)$, applying Lemma \ref{l:ind_l3}(ii), we also need the fact that $(q-\varepsilon)_3=3$ 
yields $(3,f)=1$, and so $3$ divides $|\overline G|$ only if $PGL_3^\varepsilon(q)\leq G$. Thus we may assume that $L$ satisfies the hypothesis of Lemma {\rm \ref{l:I}} and let $I$ be as in this lemma. Also observe that $\tau(L)\setminus \pi(\overline G)=\tau(L)\setminus \pi(f)$ by Lemma \ref{l:general}(i) and Lemma \ref{l:w3}. 

Let $|\overline G|$ be odd. We may assume that $G\neq L$. By Lemma \ref{l:general_ind}, if $R_i(\epsilon q)\setminus\pi(\overline G)\neq \varnothing$ for all $i\in I$, then $t(2,G)=t(2,L)$ and $\tau(G)$ is as claimed. If $R_i(\epsilon q)\setminus\pi(\overline G)=\varnothing$, then $q=8$ and $(i,\epsilon)\in\{(2,+), (1,-)\}$ by Lemma \ref{l:bz0}. By Lemma \ref{l:I}, we have $L=L_2(8)$ or $L=L_3^\varepsilon(q)$. If $L=L_2(8)$, then 
$G=\Aut L$, in which case $\tau(G)=\{7\}$ and $t(2,G)=2$. If $L=L_3(8)$, then $G=L.3$, in which case $\tau(G)=\{73\}$ and $t(2,G)=2$. If $L=U_3(8)$, then $\pi(\Out L)=\{2,3\}$ and $3\not\in \tau(L)$.  So (ii) holds.

Let Sylow $2$-subgroups of $\overline G$ be non-trivial cyclic and let $\alpha$ be some involution of $\overline G$. For every $i\in I$, we need to determine whether the set $(R_i(\epsilon q)\setminus \pi_\alpha)^*$ is empty or not. Note that the non-empty sets $R_i(\epsilon q)\setminus \pi_\alpha$ can be easily extracted from Tables \ref{tab:2even} and \ref{tab:2odd}.
 
Let $i\in I$ be such that $R_i(\epsilon q)\cap \pi_\alpha=\varnothing$. Then $(R_i(\epsilon q)\setminus \pi_\alpha)^*= R_i(\epsilon q)\setminus \pi(\overline G)$. As we showed above,  if the latter set is empty, then $L=L_2(8)$ or $L_3^\epsilon(8)$ and $i\leq 2$. Since $|\overline G|$ is even, it follows that $L=L_3^\epsilon(8)$ and $\alpha=\gamma$. But then $i=3$ by Table \ref{tab:2even}, a contradiction. So $(R_i(\epsilon q)\setminus \pi_\alpha)^*\neq \varnothing$ for such $i$.

Let $i\in I$ be such that $R_i(\epsilon q)$ neither disjoint from $\pi_\alpha$, nor is contained in it. Tables \ref{tab:2even} and \ref{tab:2odd} show this is the case if and only if $f$ is even, 
$\alpha\in\{\varphi, \varphi\gamma\}$, and $L$ is one of $L_n(q)$ with $n\geq 3$, $O_{2n}^+(q)$ with $q$ even and $n$ odd, or $E_6(q)$. Furthermore, such $i$ is unique, $i\geq 3$ and $R_i(q)\setminus \pi_\alpha=R_{2i}(q_0)$ or $R_i(q_0)$, where $q_0=p^{f/2}$, depending on whether $\alpha=\varphi$ or $\varphi\gamma$. By Lemma \ref{l:bz0}(ii), the set $R_{i}(q_0)\setminus\pi(\overline G)$ is always not empty. Suppose that $R_{2i}(q_0)\setminus\pi(\overline G)=\varnothing$. Then $i=3$ and $q_0=2$. 
It follows that $L=L_3(4)$ or $L_4(4)$. We have assumed that $L\neq L_3(4)$. If $L=L_4(4)$, then it is not hard to verify the corresponding part of Table \ref{tab:adj_even}, so we may assume that $L\neq L_4(4)$. Then $(R_i(\epsilon q)\setminus \pi_\alpha)^*\neq \varnothing$ for such $i$ either.

If $p=2$, then $\alpha$ is conjugate to $\varphi$, or $\varphi\gamma$, or $\gamma$ and looking through Table \ref{tab:2even}, we verify Table \ref{tab:adj_even} and the corresponding part of Table \ref{tab:adj_ex}.

In odd characteristic, we give a detailed proof for the case when $L=L_n^\varepsilon(q)$ and $n$ is even, the other cases are similar. 

Let $L=L_n^\varepsilon(q)$, where $n\geq 4$ is even and $q$ is odd. 
Suppose that $|(\Outdiag L\rtimes \langle \gamma\rangle)\cap\overline G|$ is even. If $\delta\in\overline G$, then $t(2,G)=1$ by Lemma \ref{l:general}(iii). If $\delta\not\in\overline G$, then $\alpha$ is conjugate to $\gamma$ or $\gamma'$. If $(n)_2\geq (q-\varepsilon)_2$, then the only non-empty set $R_i(\varepsilon q)\setminus \pi_\alpha$ is $R_{n-1}(\varepsilon q)$, and $R_{n-1}(\varepsilon q)^*\neq \varnothing$ by the above argument. So in this case $\tau(G)=R_{n-1}(\varepsilon q)^*$ and $t(2,G)=2$. If $(n)_2<(q-\varepsilon)_2$, then $t(2,G)=1$ unless $\varepsilon=-$, $(n)_2=2$ and $\alpha$ is conjugate to $\gamma'$, in which case the only non-empty set $R_i(\varepsilon q)\setminus \pi_\alpha$ is $R_n(-q)$.

Let  $|(\Outdiag L\rtimes \langle\gamma\rangle)\cap\overline G|$ be odd. Then $\varepsilon=+$ and $f$ is even.  If $(n)_2\geq (q-1)_2$, then $\overline G$ contains a conjugate of $\varphi$ or $\varphi\gamma$ by Lemma \ref{l:involutions}. Taking $\tau(L)\setminus\pi_\alpha$ for $\alpha=\varphi, \varphi\gamma$ from Table~\ref{tab:2odd}, we get the corresponding part of Table \ref{tab:adj_odd}. If $(n)_2<(q-1)_2$, then Table~\ref{tab:2odd} states that $\tau(L)\setminus\pi_\alpha=R_n(q)$, and this gives the  corresponding part of Table \ref{tab:adj_odd}.

The proof is complete.
\end{proof}

\section{Almost simple groups with alternating or sporadic socle}

In this section, for completeness, we give the solutions to Problems \ref{p:asWr'} and \ref{p:ast2} for almost simple groups with alternating or sporadic socle. 

Recall that $\Aut A_n = S_n$ if $n=5$ or $n\geq 7$. Also  $A_6\simeq L_2(9)$ and, dealing with almost simple groups, it is more convenient to regard $A_6$ as a group of Lie type. Denote by $\mathbb P$ the set of all primes.

\begin{prop}\label{prop:alt} Let $n\geq 5$ and $r$ be an odd prime. Then 
\begin{enumerate}
 \item $A_n$ is a $W_r$-group if and only if $n\in\{r,r+1,r+2,r+3\}$ and $S_n$ is a $W_r$-group if and only if $n\in\{r,r+1\}$. 
\item $\tau(A_n)=\{n, n-1,n-2,n-3\}\cap \mathbb P$. In particular, $t(2,A_n)=3$ if there are two primes in the set $\{n, n-1,n-2,n-3\}$, $t(2,A_n)=2$ if there is one prime in the  set $\{n, n-1,n-2,n-3\}$, and $t(2,A_n)=1$ otherwise.
\item $\tau(S_n)=\{n, n-1\}\cap \mathbb P$. In particular, $t(2,S_n)=2$ if there is a prime in the set $\{n, n-1\}$ and $t(2,S_n)=1$ otherwise.
\end{enumerate}
\end{prop}

\begin{proof} This is obvious since an odd $r\in\pi(A_n)$ is adjacent to $2$ in $GK(A_n)$ if and only if $4+r\leq n$. For $S_n$, the corresponding inequality is $2+r\leq n$. 
\end{proof}

\begin{prop}\label{prop:spor} Let $G$ be an almost simple  group with sporadic socle. Then $\tau(G)$ and $t(2,G)=|\tau(G)|+1$ are as in Table {\rm \ref{tab:wr_spor}}. 
 \end{prop}
\begin{proof} See \cite{85Atlas}.
\end{proof}

\begin{table}[h]

\caption{The set $\tau(G)$ for almost simple sporadic groups}\label{tab:wr_spor}
\begin{tabular}{|l|l|l|}
\hline 
$G$&$\tau(G)$&$t(2,G)$\\
\hline
{$M_{11}\simeq \Aut M_{11}$} & {$5$, $11$}&$3$\\
{$M_{12}$, $\Aut M_{12}$} & {$11$}&$2$\\
{$M_{22}$} & {$5$, $7$, $11$}&$4$ \\
{$\Aut M_{22}$} & {$11$}&$2$\\
{$M_{23}\simeq \Aut M_{23}$} & {$5$, $11$, $23$}&$4$\\
{$M_{24}\simeq \Aut M_{24}$} & {$11$, $23$}&$3$\\
{$J_1\simeq \Aut J_1$} & {$7$, $11$, $19$}&$4$\\
{$J_2$} & {$7$}&$2$\\
{$\Aut J_2$} & $\varnothing$&$1$\\
{$J_3$} & {$17$, $19$}&$3$\\
{$\Aut J_3$} & {$19$}&$2$ \\
{$J_4\simeq \Aut J_4$} & {$23$, $29$, $31$, $37$, $43$} &$6$\\
{$HS$} & {$7$, $11$}&$3$ \\ 
{$\Aut HS$} & {$11$}&$2$\\
{$McL$} & {$11$}&$2$\\
{$\Aut McL$)} & $\varnothing$&$1$\\
{$He$, $\Aut He$} & {$17$}&$2$ \\

\hline
\end{tabular}
\begin{tabular}{|l|l|l|}
\hline 
$G$&$\tau(G)$&$t(2,G)$\\
\hline
{$Ru\simeq \Aut Ru$} & {$29$}&$2$ \\
{$Suz$} & {$11$, $13$}&$3$ \\
{$\Aut Suz$} & {$13$}&$2$ \\  
{$O'N$} & {$11$, $19$, $31$} &$4$\\
{$\Aut O'N$}& {$31$}&$2$\\
{$Co_1\simeq \Aut Co_1$} & {$23$}&$2$\\
{$Co_2\simeq \Aut Co_2$} & {$11$, $23$}&$3$\\
{$Co_3\simeq \Aut Co_3$} & {$23$}&$2$ \\
{$Fi_{22}$, $\Aut Fi_{22}$} & {$13$}&$2$\\
{$Fi_{23}\simeq \Aut Fi_{23}$} & {$17$, $23$}&$3$\\
{$Fi_{24}'$} & {$17$, $23$, $29$}&$4$\\
{$\Aut Fi_{24}'$} & {$29$}&$2$\\
{$HN$, $\Aut HN$} & {$19$}&$2$\\
{$Ly\simeq \Aut Ly$} & {$31$, $37$, $67$}&$4$ \\
{$Th\simeq \Aut Th$} & {$13$, $19$, $31$}&$4$\\
{$B\simeq \Aut B$} & {$31$, $47$} &$3$\\
{$M \simeq \Aut M$} & {$29$, $41$, $59$ $71$}&$5$ \\
\hline
\end{tabular}
\end{table}

\section*{Acknowledgments}

The authors are grateful to Natalia Maslova for helpful discussions on this work. 
The work of M. A. Grechkoseeva was carried out within the framework of the Sobolev
Institute of Mathematics state contract (project FWNF-2026-0017). Yang Nanying is supported by Fundamental Research Funds for the Central Universities-President’s Fund (JUSRP202406006).


\end{document}